\documentclass[a4paper,11pt]{amsart}
\usepackage[T1]{fontenc}
\usepackage{lmodern}
\usepackage{microtype}
\usepackage{amsmath,amssymb,amsthm,mathtools}
\usepackage{booktabs}
\usepackage{array}
\usepackage{enumitem}
\usepackage{tikz}
\usetikzlibrary{arrows.meta,decorations.pathreplacing}
\usepackage{placeins}
\usepackage{needspace}
\usepackage[hidelinks]{hyperref}
\hypersetup{
  pdftitle={Semiprime Rigidity and Growth along Powers in Integer-Sided Cyclic Polygons},
  pdfauthor={Felix Huber},
  pdfsubject={Semiprime rigidity, growth along prime powers, and parity of integer-sided cyclic polygons},
  pdfkeywords={integer-sided cyclic polygons, semiprime radii, Gaussian integers, Eisenstein integers, Farey neighbors, parity}
}

\newtheorem{mainthm}{Theorem}

\newtheorem{theorem}{Theorem}[section]
\newtheorem{proposition}[theorem]{Proposition}
\newtheorem{lemma}[theorem]{Lemma}
\newtheorem{corollary}[theorem]{Corollary}
\theoremstyle{definition}
\newtheorem{definition}[theorem]{Definition}
\newtheorem{example}[theorem]{Example}
\theoremstyle{remark}
\newtheorem{remark}[theorem]{Remark}

\newcommand{\Z}{\mathbb Z}
\newcommand{\Q}{\mathbb Q}
\newcommand{\OK}{\mathcal O_K}
\newcommand{\Norm}{\operatorname{N}}

\title[Semiprime rigidity and power growth]
{Semiprime Rigidity and Growth along Powers\\in Integer-Sided Cyclic Polygons}
\author{Felix Huber}
\address{Independent Researcher, Lucerne, Switzerland}
\email{felix.68@gmx.ch}
\thanks{ORCID: 0009-0005-1568-1579.}
\date{November 2026}

\begin{document}

\begin{abstract}
For a positive integer $n$, let $a(n)$ be the largest number of sides of a
simple cyclic polygon of circumradius $n$ whose side lengths are positive
integers. For primes $p,q$ we show that $a(pq)\in\{6,8,18\}$ and characterize
each value by conditions on Gaussian and Eisenstein integers. For every
integer $N\ge2$, on the other hand, $a(N^e)\to\infty$ and
$\limsup_{e\to\infty}a(N^e)/e\ge2\sqrt2$; a quadratic chord tower combined with
Farey neighbors gives explicit polygons, for instance $a(8192)\ge38$. We also
show that an odd-sided polygon of radius $n$ exists if and only if $n$ has a
prime divisor congruent to $1\pmod4$, and that $a(n)$ is even whenever at most
one distinct prime divisor of $n$ is congruent to $1\pmod4$;
it is also even when the Gaussian part of $n$ is the product
of two distinct primes congruent to $1\pmod4$.
Moreover, $a(N^e)=O_N(e^{\omega(N)})$ for every fixed $N\ge2$, and
$a(q^e)=\Theta(e)$ for primes $q\equiv1\pmod4$; exceptional semiprime
extremisers have unique side-length multisets. The proofs use a
splitting lemma for angle relations by squarefree class, arithmetic in
imaginary quadratic fields, bounded zero-sum problems, and a one-dimensional
Farey argument. No computer-assisted step is used.
\end{abstract}

\maketitle

\section{Introduction}\label{sec:intro}

Fix a positive integer $n$ and consider simple polygons inscribed in a circle
of radius $n$ all of whose sides have positive integer length. How many sides
can such a polygon have? A regular hexagon, with side $n$, always gives six.
At radius $65$ there is an octagon with sides $32,50,50,66,32,50,50,66$
(Example~\ref{ex:gauss}); at radius $91$ there is an eighteen-gon whose sides
$7,26,61$ each occur six times (Example~\ref{ex:eis}); and at radius
$2^{13}=8192$ the construction of Section~\ref{sec:powers} gives a polygon
with $38$ sides.

A chord of length $s$ in a circle of radius $n$ subtends a central angle whose
half is
\begin{equation}\label{eq:alpha}
 \alpha_n(s)=\arcsin\!\left(\frac{s}{2n}\right)\in(0,\pi/2]\qquad(1\le s\le2n),
\end{equation}
and a polygon closes up precisely when such angles satisfy one exact additive
relation (Lemma~\ref{lem:geometry}). Which angles occur, and which additive
relations hold among them, depends on the arithmetic of $n$. Much of the
classical literature on cyclic polygons concerns areas and rationality
conditions~\cite{Robbins,MaleyRobbinsRoskies,BuchholzMacDougall}, and
trigonometric Diophantine relations have a separate
tradition~\cite{ConwayJones,CRS}. Here we study the \emph{maximum number of
integer chords that close up} at a prescribed integer circumradius.

\begin{definition}\label{def:a}
Let $a(n)$ be the greatest integer $k\ge3$ for which there is a simple
cyclic polygon with $k$ distinct vertices, circumradius $n$, and positive
integer side lengths.
\end{definition}
The regular hexagon has integer side length $n$, so $a(n)\ge6$.
The maximum exists: every side subtends a central angle of at least
$2\alpha_n(1)$, so the full turn bounds the number of sides.
The sequence is OEIS A398490; its associated counting triangle is
A398489~\cite{OEIS}.

\subsection*{Results}
\textbf{Semiprime radii.} Our first result determines $a(n)$ when $n$ is a
product of two primes, not necessarily distinct. The exceptional values are
governed by angles arising from the Eisenstein and Gaussian integers. To
state the criterion, let $r\equiv1\pmod3$ be prime and let $(s_r,u_r)$ be the
solution of
\begin{equation}\label{eq:eisenstein-representation}
 s_r^2+3u_r^2=4r^2,\qquad 0<s_r<r,\quad u_r>0;
\end{equation}
existence and uniqueness are proved in Lemma~\ref{lem:eisdata}. Write
\begin{equation}\label{eq:xr}
 x_r=\arcsin\!\left(\frac{s_r}{2r}\right)\in(0,\pi/6).
\end{equation}

\begin{mainthm}[Semiprime radii]\label{thm:main}
Let $p\le q$ be primes and consider the conditions
\begin{itemize}[nosep,leftmargin=2.6em]
\item[(E)] $p\equiv q\equiv1\pmod3$ (so $p,q$ are odd) and $x_p+x_q<\pi/6$;
\item[(G)] $p<q$ and $p\equiv q\equiv1\pmod4$.
\end{itemize}
Then
\[
 a(pq)=
 \begin{cases}
 18,&\text{if (E) holds},\\[1mm]
 8,&\text{if (G) holds but (E) does not},\\[1mm]
 6,&\text{otherwise.}
 \end{cases}
\]
Moreover, $x_p+x_q<\pi/6$ is equivalent to the integer inequality
\begin{equation}\label{eq:integer-criterion}
 3(u_ps_q+u_qs_p)<3u_pu_q-s_ps_q.
\end{equation}
In particular,
\[
 a(pq)\in\{6,8,18\},\qquad a(2p)=a(3p)=6
\]
for every prime $p$.
\end{mainthm}

For a semiprime radius, each integer chord determines a vector of
prime-ideal valuations in a small box, and an exact closing relation becomes a
zero-sum relation among these vectors. In an imaginary quadratic field with
only the units $\pm1$ this does not allow more than six sides; the additional
roots of unity in the Gaussian and Eisenstein fields allow eight and eighteen
sides. For example, $a(49)=18$, $a(65)=8$, $a(91)=18$, while $a(361)=6$
(Section~\ref{sec:examples}). For $n=1$ and every prime $n$ one has $a(n)=6$
(Remark~\ref{rem:primes}). Semiprime radii are therefore the first level,
classified by the number of prime factors counted with multiplicity, at which
the maximum can exceed six.

\medskip
\textbf{Powers of a fixed integer.} Restricting the prime divisors of the
radius does not bound $a(n)$: successive powers of a fixed integer admit
polygons with arbitrarily many sides.

\begin{mainthm}[Growth along powers]\label{thm:power-growth}
For every integer $N\ge2$,
\[
 a(N^e)\longrightarrow\infty\qquad(e\to\infty),
\]
and more quantitatively
\begin{equation}\label{eq:limsup-growth}
 \limsup_{e\to\infty}\frac{a(N^e)}{e}\ge2\sqrt2.
\end{equation}
\end{mainthm}

The proof builds integer chords from powers of one imaginary quadratic
integer; neighboring fractions on a Stern--Brocot path give the exact
cancellation of angles that closes the polygon. A general Steinitz bound
also gives $a(N^e)=O_N(e^{\omega(N)})$, and for primes $q\equiv1\pmod4$
we obtain the sharp order $a(q^e)=\Theta(e)$. For $N=2$ it yields a
$38$-gon at radius $8192=2^{13}$ (Example~\ref{ex:8192}).

\medskip
\textbf{Parity.} Theorem~\ref{thm:odd-polygon} shows that an integer-sided
cyclic polygon with an odd number of sides exists at radius $n$ if and only if
$n$ has a prime divisor congruent to $1$ modulo $4$. Whether an
\emph{extremal} polygon can have an odd number of sides is a subtler question.
Define
\[
\omega_1(n)=\#\{q:q\mid n,\ q\text{ prime},\ q\equiv1\pmod4\},
\]
where distinct primes are counted only once.

\begin{mainthm}[Parity in Gaussian rank one]\label{thm:parity}
If $\omega_1(n)\le1$, then $a(n)$ is even. In particular,
\[
 a(p^e)\equiv0\pmod2
\]
for every prime $p$ and every $e\ge1$.
\end{mainthm}

A further consequence (Corollary~\ref{cor:gaussian-part-parity})
extends this evenness conclusion to radii whose Gaussian part
is the product of two distinct primes congruent to $1\pmod4$,
regardless of additional factors $2$ or primes congruent to $3\pmod4$.

\subsection*{Organization}
All three theorems rest on exact additive relations among the angles
$\alpha_n(s)$. Section~\ref{sec:geometry} gives the closing criterion.
Section~\ref{sec:splitting} shows that such relations split by squarefree
class, which also yields the odd-sided criterion.
Section~\ref{sec:semiprime} proves Theorem~\ref{thm:main}: zero-sum relations
bound the number of chords per unit angle in each class, and explicit integer
constructions attain the bounds in the Gaussian and Eisenstein classes.
Sections~\ref{sec:powers} and~\ref{sec:parity} prove
Theorems~\ref{thm:power-growth} and~\ref{thm:parity}.
Section~\ref{sec:examples} contains examples and Section~\ref{sec:outlook}
open questions.

\section{The geometric criterion}\label{sec:geometry}

\begin{lemma}[Minor- and major-arc criterion]\label{lem:geometry}
A multiset of $k\ge3$ positive integers $\{s_1,\ldots,s_k\}$, with $s_i\le2n$, is
realizable as the side-length multiset of a simple cyclic polygon of radius
$n$ if and only if either
\begin{equation}\label{eq:minor}
 \sum_{i=1}^k\alpha_n(s_i)=\pi,
\end{equation}
or
\begin{equation}\label{eq:major}
 \alpha_n(s_j)=\sum_{i\ne j}\alpha_n(s_i)
\end{equation}
for some $j$.
\end{lemma}

\begin{proof}
For points on a circle, the boundary edges of a simple polygon follow
their circular order, up to reversal. Indeed, if an edge joined two
nonconsecutive vertices, the remaining vertices on the two open arcs
between its endpoints would lie on opposite sides of that chord.
The remaining boundary path starts and ends at the chord's endpoints.
Its interior visits vertices on both open arcs, so one of its edges must
join vertices on opposite sides of the chord and therefore cross it;
this contradicts simplicity. Thus the successive boundary arcs form
an ordered partition of the full circle and their central angles
sum to $2\pi$.
If every such arc has length at most $\pi$, its central angle is
$2\alpha_n(s_i)$, giving \eqref{eq:minor}. Otherwise exactly one
arc has length greater than $\pi$ (two such arcs would already exceed $2\pi$),
and its central angle is $2\pi-2\alpha_n(s_j)$; the other arcs are minor,
and the central angles sum to $2\pi$ exactly when \eqref{eq:major} holds.

Conversely, suppose $k\ge3$ and one of these equations holds.
Arrange on the circle consecutive arcs with central angles $2\alpha_n(s_i)$
in the first case, and $2\alpha_n(s_i)$ $(i\ne j)$ together with
$2\pi-2\alpha_n(s_j)>0$ in the second. Each angle is strictly positive and
their sum is $2\pi$; therefore the cumulative endpoints give $k$ distinct
circle vertices in cyclic order. Joining consecutive vertices yields a
convex, hence simple, cyclic $k$-gon, and the chord subtending an arc of
central angle $\theta\in(0,2\pi)$ has length $2n\sin(\theta/2)$, which equals
the prescribed $s_i$ in both cases.
\end{proof}

\begin{figure}[htbp]
\centering
\begin{minipage}[t]{0.48\textwidth}
\centering
\begin{tikzpicture}
\def\rad{1.5}
\coordinate (O) at (0,0);
\draw[gray!60] (O) circle (\rad);
\foreach \a/\v in {40/P,140/Q,215/R,290/S,355/T}{\coordinate (\v) at (\a:\rad);}
\draw[thick] (P)--(Q)--(R)--(S)--(T)--cycle;
\draw[gray,densely dashed] (O)--(P) (O)--(Q);
\draw[line width=1.6pt] (P)--(Q);
\draw (40:0.22) arc (40:140:0.22);
\node[font=\scriptsize,fill=white,inner sep=0.8pt] at (90:0.52) {$2\alpha_n(s_i)$};
\foreach \v in {P,Q,R,S,T}{\fill (\v) circle (1.2pt);}
\fill (O) circle (1pt); \node[font=\scriptsize,below=1pt] at (O) {$O$};
\node[font=\scriptsize,above right] at (P) {$P$};
\node[font=\scriptsize,above left] at (Q) {$Q$};
\path (-2.0,-2.1) rectangle (2.0,2.0);
\end{tikzpicture}

\smallskip
{\small (a) All boundary arcs are minor.}
\end{minipage}\hfill
\begin{minipage}[t]{0.48\textwidth}
\centering
\begin{tikzpicture}
\def\rad{1.5}
\coordinate (O) at (0,0);
\foreach \a/\v in {25/P,70/B,110/C,140/D,155/Q}{\coordinate (\v) at (\a:\rad);}
\draw[gray!60] (O) circle (\rad);
\draw[thick] (P)--(B)--(C)--(D)--(Q)--cycle;
\draw[gray,densely dashed] (O)--(P) (O)--(Q);
\draw[line width=1.6pt] (P)--(Q);
\draw (25:0.22) arc (25:155:0.22);
\node[font=\scriptsize,fill=white,inner sep=0.8pt] at (90:0.40) {$2\alpha_n(s_j)$};
\draw[line width=1.5pt,dash pattern=on 4pt off 2pt] (155:\rad) arc (155:385:\rad);
\node[font=\scriptsize] at (270:\rad+0.32) {$2\pi-2\alpha_n(s_j)$};
\foreach \v in {P,B,C,D,Q}{\fill (\v) circle (1.2pt);}
\fill (O) circle (1pt); \node[font=\scriptsize,below=1pt] at (O) {$O$};
\node[font=\scriptsize,above right] at (P) {$P$};
\node[font=\scriptsize,above left] at (Q) {$Q$};
\path (-2.0,-2.1) rectangle (2.0,2.0);
\end{tikzpicture}

\smallskip
{\small (b) One boundary arc is major.}
\end{minipage}
\caption{Minor and major configurations in Lemma~\ref{lem:geometry}. The heavy
chord $PQ$ has length $s_i$ in (a) and $s_j$ in (b); dashed segments join $P$
and $Q$ to the center $O$. In (b) the boundary arc from $Q$ to $P$ (dashed
circular arc) has central angle $2\pi-2\alpha_n(s_j)$, while the other sides
together subtend the central angle $2\alpha_n(s_j)$ of $PQ$.
The diagram was drafted with AI-assisted TikZ code from this geometric construction.}
\label{fig:arc-geometry}
\end{figure}
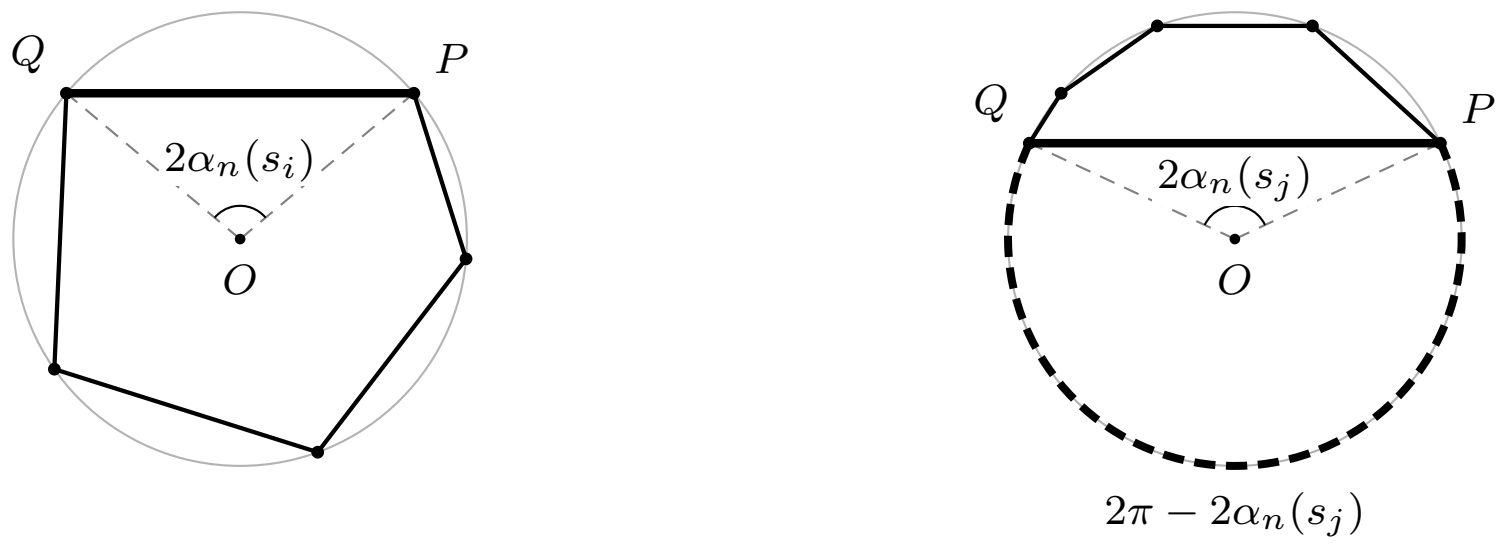

\begin{remark}\label{rem:diameter}
A side of length $2n$ (a \emph{diameter}) has $\alpha_n(2n)=\pi/2$. At most one
side can be a diameter, because the other $k-1\ge2$ half-angles are positive.
A relation \eqref{eq:major} with $\alpha_n(s_j)=\pi/2$ is the same as
\eqref{eq:minor}. Hence every polygon is of exactly one of three kinds:
a polygon \emph{with a diameter}, a \emph{minor-arc polygon without diameter}
(relation \eqref{eq:minor}), or a \emph{major-arc polygon} (relation
\eqref{eq:major} with $\alpha_n(s_j)<\pi/2$). Figure~\ref{fig:arc-geometry}
illustrates the last two kinds.
\end{remark}

\begin{definition}[Blocks]\label{def:block}
Put $H=\pi/2$. A \emph{block} at radius $n$ is a finite nonempty multiset of
half-angles $\alpha_n(s)$ of \emph{non-diameter} sides ($1\le s<2n$); its
\emph{total} is the sum of its members. A block is a \emph{$\pi$-block} or an
\emph{$H$-block} if its total is $\pi$ or $H$, respectively, and a
\emph{rational block} if its total lies in $\pi\mathbb Q$. A block is
\emph{pure} in a squarefree class $d$ (defined in Section~\ref{sec:splitting})
if all its members lie in class $d$.
\end{definition}
Thus \eqref{eq:minor} says that the sides of a minor-arc polygon without
diameter form a $\pi$-block, and a polygon with a diameter consists of the
diameter together with an $H$-block (Remark~\ref{rem:diameter}); relation
\eqref{eq:major} is the \emph{major-arc relation}.

We shall also use the strict monotonicity
\begin{equation}\label{eq:monotonicity}
 \frac{d}{ds}\alpha_n(s)=\frac1{\sqrt{4n^2-s^2}}>0.
\end{equation}
Thus a half-angle uniquely determines its side length at fixed radius.

\begin{proposition}[Scaling]\label{prop:scaling}
For all positive integers $m,n$,
\[
 a(mn)\ge a(n).
\]
\end{proposition}

\begin{proof}
Scale every side length of a polygon of radius $n$ by $m$. All central
angles are unchanged, and the new polygon has radius $mn$.
\end{proof}

\section{Squarefree classes and elementary splitting}\label{sec:splitting}
We group chords by the squarefree part of the complementary chord equation
and show that rational angle relations split across these classes. This is the
arithmetic input for Theorems~\ref{thm:main} and~\ref{thm:parity}.

For a non-diameter side $s$ at radius $n$, write
\begin{equation}\label{eq:squarefree-class}
 4n^2-s^2=du^2,
\end{equation}
where $d>0$ is squarefree and $u$ is positive; we say that $s$ (or its
half-angle) lies in the \emph{class} $d$. Recall $H=\pi/2$
and set $\beta_n(s)=H-\alpha_n(s)\in(0,H)$. Then
\begin{equation}\label{eq:tanbeta}
 \tan\beta_n(s)=\frac{u}{s}\sqrt d.
\end{equation}
In particular, $\tan^2\beta_n(s)=du^2/s^2$ is rational, and the tangent
belongs to a fixed squarefree class $d$. Angles with rational squared tangent
are the subject of Conway, Radin, and Sadun~\cite{CRS}; we prove the splitting
property needed below directly, using elementary ideal theory and Galois
characters. The splitting assertion below is a special case of the Splitting Theorem
of Conway--Radin--Sadun~\cite[Theorem~3]{CRS}; we retain our
self-contained character-orthogonality proof.
The classes $d=1$ (\emph{Gaussian}) and $d=3$
(\emph{Eisenstein}) are special because $\mathbb Q(\sqrt{-d})$ then has more
than two roots of unity; all other classes are called \emph{ordinary}.

For later normalization, put
\begin{equation}\label{eq:primitive-norm}
 \kappa=\gcd(s,u),\qquad A=\frac{s}{\kappa},\qquad B=\frac{u}{\kappa},
 \qquad M=\frac{2n}{\kappa}.
\end{equation}
Because $\kappa^2$ divides $s^2+du^2=(2n)^2$, one has $\kappa\mid2n$;
thus $M$ is a positive integer and $M\mid2n$. Moreover,
\begin{equation}\label{eq:norm-equation}
 \gcd(A,B)=1,
 \qquad A^2+dB^2=M^2.
\end{equation}

\subsection{Ideal components}
With the normalization \eqref{eq:primitive-norm} and
\eqref{eq:norm-equation}, put $K=\Q(\sqrt{-d})$,
$z=A+B\sqrt{-d}$, and $w=z/M$. Then $|w|=1$, $\arg w=\beta_n(s)=H-\alpha_n(s)\in(0,H)$, and
\begin{equation}\label{eq:idealnorm}(z)(\bar z)=(M)^2.\end{equation}
For each rational prime $r$ that splits as $(r)=\mathfrak r\bar{\mathfrak r}$, fix once and for all an ordering of the two prime ideals and put
\[
 v_r(w)=\tfrac12\bigl(\operatorname{ord}_{\mathfrak r}(w)-\operatorname{ord}_{\bar{\mathfrak r}}(w)\bigr)
=\tfrac12\bigl(\operatorname{ord}_{\mathfrak r}(z)-\operatorname{ord}_{\bar{\mathfrak r}}(z)\bigr).
\]
The equality \eqref{eq:idealnorm} makes the two exponents of $z$ sum to $2\operatorname{ord}_rM$; hence their difference is even and $v_r(w)\in\Z$. The fractional ideal $(w)$ therefore has exponents $v_r(w),-v_r(w)$ at $\mathfrak r,\bar{\mathfrak r}$. At inert and ramified primes the exponents of $(w)$ vanish, and so do those at primes not dividing $M$. In particular, only split primes dividing $2n$ may contribute. We call $v_r(w)$ the \emph{$r$-component} of the chord. The component vectors of different chords are always computed with the same fixed orientation of each split prime-ideal pair.

\begin{lemma}[Component criterion]\label{lem:component}
Let $w_1,\dots,w_m$ be the numbers attached as above to non-diameter sides
(possibly at different radii) of one \emph{fixed class} $d$, with half-angles
$\alpha_i$, and let $\varepsilon_i\in\mathbb Z$. Then
\begin{equation}\label{eq:zero}
 \sum_i\varepsilon_i\alpha_i\in\pi\Q
 \quad\Longleftrightarrow\quad
 \sum_i\varepsilon_i v_r(w_i)=0\text{ for every split prime }r.
\end{equation}
Moreover, in that case $\prod_iw_i^{\varepsilon_i}$ is a root of unity of $K$.
The unit group of $\OK$ is $\{\pm1\}$ for $d\ne1,3$, $\mu_4$ for $d=1$, and $\mu_6$ for $d=3$.
\end{lemma}

\begin{proof}
The product $W=\prod_iw_i^{\varepsilon_i}\in K^\times$ has absolute value $1$ and
argument $\sum_i\varepsilon_i\beta_i=(\sum_i\varepsilon_i)H-\sum_i\varepsilon_i\alpha_i$.
Since $H\in\pi\Q$, the left-hand condition in \eqref{eq:zero} says exactly that $\arg W\in\pi\Q$, i.e.\ that
$W$ is a root of unity. The fractional ideal of $W$ is
$\prod_r(\mathfrak r/\bar{\mathfrak r})^{\sum_i\varepsilon_iv_r(w_i)}$ by the preceding paragraph, and this vanishes
if and only if all exponents vanish. If $W$ is a root of unity its ideal is trivial, which gives the forward implication. Conversely, if the ideal of
$W$ is trivial, then $W$ is a unit of $\OK$, hence a root of unity because the unit group is finite.
\end{proof}

\begin{lemma}[Elementary splitting of squarefree classes]\label{lem:splitting}
Let $\alpha_i$ be finitely many non-diameter integer-chord half-angles,
possibly at different radii, and let $\varepsilon_i\in\mathbb Z$.
If $\sum_i\varepsilon_i\alpha_i\in\pi\mathbb Q$, then for every
squarefree $d$ the restricted sum
$\sum_{i:d_i=d}\varepsilon_i\alpha_i$ belongs to $\pi\mathbb Q$.
\end{lemma}
\begin{proof}
As $\beta_i=H-\alpha_i$, the hypothesis also gives
$\sum_i\varepsilon_i\beta_i\in\pi\mathbb Q$.
For every occurring squarefree $d$, put
$K_d=\mathbb Q(\sqrt{-d})$ and
$W_d=\prod_{i:d_i=d}w_i^{\varepsilon_i}\in K_d$.
We have $|W_d|=1$ and
$\prod_d W_d=\zeta$, a root of unity.
Let $L$ be the compositum of the fields $K_d$ and
$G=\operatorname{Gal}(L/\mathbb Q)$. It is a finite elementary
abelian $2$-group. Each field $K_d$ defines a character
$\chi_d:G\to\{\pm1\}$: a Galois automorphism acts on
$K_d$ as the identity when $\chi_d=1$ and as complex conjugation when
$\chi_d=-1$. Distinct squarefree $d$ give distinct characters, since
$\ker\chi_d=\operatorname{Gal}(L/K_d)$ determines $K_d$.
Since $\overline{W_d}=W_d^{-1}$,
$\sigma(W_d)=W_d^{\chi_d(\sigma)}$.
Fix a class $d_0$. Apply every $\sigma\in G$ to
$\prod_dW_d=\zeta$, take the $\chi_{d_0}(\sigma)$-th power,
and multiply over $\sigma$; this yields
\[
\prod_{\sigma\in G}\sigma(\zeta)^{\chi_{d_0}(\sigma)}
 =\prod_d W_d^{\sum_{\sigma\in G}\chi_{d_0}(\sigma)\chi_d(\sigma)}
 =W_{d_0}^{|G|},
\]
by the orthogonality of distinct characters of the finite abelian
group $G$. The left side is a root of unity; hence $W_{d_0}$ is a
root of unity. Thus
$\sum_{i:d_i=d_0}\varepsilon_i\beta_i\in\pi\mathbb Q$,
and subtracting from the corresponding integral multiple of $H$
proves the assertion for the $\alpha_i$.
\end{proof}

The next lemma locates the rational sub-sums of Lemma~\ref{lem:splitting}. It
is used for the parity principle, for major-arc relations, and in
Section~\ref{sec:parity}.

\begin{lemma}[Unit lattices of rational sub-sums]\label{lem:unitlattice}
Let $\alpha_i$ and $\varepsilon_i$ be as in Lemma~\ref{lem:splitting}, with
$\sum_i\varepsilon_i\alpha_i\in\pi\mathbb Q$. For a squarefree class $d$ put
\[
 \mathcal A_d=\sum_{i:\,d_i=d}\varepsilon_i\alpha_i,\qquad
 m_d=\sum_{i:\,d_i=d}\varepsilon_i,\qquad U=\frac\pi6 .
\]
Then $\mathcal A_d\in\pi\Q$ and:
\begin{enumerate}[nosep,label=(\roman*)]
\item if $d\ne1,3$, then $\mathcal A_d\in H\Z$ and $\mathcal A_d/H\equiv m_d\pmod2$;
\item if $d=1$, then $\mathcal A_d\in H\Z$;
\item if $d=3$, then $\mathcal A_d\in U\Z$ and $\mathcal A_d/U\equiv m_d\pmod2$.
\end{enumerate}
\end{lemma}
\begin{proof}
By Lemma~\ref{lem:splitting}, $\mathcal A_d\in\pi\Q$, so
$\mathcal B_d=\sum_{i:d_i=d}\varepsilon_i\beta_i=m_dH-\mathcal A_d\in\pi\Q$, and by
Lemma~\ref{lem:component} the number $W_d=\prod_{d_i=d}w_i^{\varepsilon_i}$ is a root of
unity of $\Q(\sqrt{-d})$ with argument $\mathcal B_d$ (modulo $2\pi$). Hence
$\mathcal B_d\in\pi\Z=2H\Z$ for $d\ne1,3$, $\mathcal B_d\in H\Z$ for $d=1$, and
$\mathcal B_d\in(\pi/3)\Z=2U\Z$ for $d=3$. Now use $\mathcal A_d=m_dH-\mathcal B_d$,
and for $d=3$ use $H=3U$ and $3m_d\equiv m_d\pmod 2$.
\end{proof}

The following parity principle gives the odd-sided criterion and isolates the
Gaussian class as the only obstruction to evenness.

\begin{lemma}[Non-Gaussian parity principle]\label{lem:nongaussian-parity}
Let $\alpha_1,\ldots,\alpha_m$ be half-angles of non-diameter integer chords
of the same radius, none belonging to the Gaussian class $d=1$. Let
$\varepsilon_i\in\{1,-1\}$, and suppose
\[
 \sum_{i=1}^m \varepsilon_i\alpha_i=tH,
 \qquad H=\frac\pi2,
 \qquad t\in\Z.
\]
Then
\[
 \sum_{i=1}^m\varepsilon_i\equiv t\pmod2.
\]
\end{lemma}

\begin{proof}
Lemma~\ref{lem:unitlattice} applies with the sums $\mathcal A_d$, $d\ne 1$. For
$d\ne1,3$ we have $\mathcal A_d/H\equiv m_d\pmod2$. For $d=3$ we have only
$\mathcal A_3\in U\Z$; but $\mathcal A_3=tH-\sum_{d\ne3}\mathcal A_d\in H\Z$, so
$\mathcal A_3=\ell H=3\ell U$ for an integer $\ell$, and
$\ell\equiv3\ell=\mathcal A_3/U\equiv m_3\pmod2$. Thus
$\mathcal A_d/H\equiv m_d\pmod2$ for every occurring class $d$. Summing over
$d$ gives $t=\sum_d\mathcal A_d/H\equiv\sum_dm_d=\sum_i\varepsilon_i\pmod2$.
\end{proof}

\begin{proposition}[Parity without Gaussian chords]\label{cor:nongaussian-polygons}
At any radius $n$, a simple integer-sided cyclic polygon containing no
non-diameter side of squarefree class $d=1$ has an even number of sides.
\end{proposition}
\begin{proof}
Let $k$ be the number of sides. By Remark~\ref{rem:diameter} there are three kinds of polygon.
\emph{Minor-arc polygon without diameter:}
Lemma~\ref{lem:nongaussian-parity} with all signs $+1$ and $t=2$ gives
$k\equiv2\equiv0\pmod2$.
\emph{Major-arc polygon:} give the target $\alpha_j$ the sign $+1$ and each other
side the sign $-1$; the signed sum is $0$, so $t=0$ and
$1-(k-1)=2-k\equiv0$, i.e.\ $k$ is even.
\emph{Polygon with a diameter:} the remaining $k-1$ non-diameter half-angles
sum to $H$, so Lemma~\ref{lem:nongaussian-parity} with all signs $+1$ and $t=1$
shows that $k-1$ is odd, hence $k$ is even.
\end{proof}

\begin{theorem}[Odd-sided polygons]\label{thm:odd-polygon}
For a positive integer $n$, the following are equivalent:
\begin{enumerate}[label=(\roman*)]
\item there exists an odd-sided integer-sided cyclic polygon of radius $n$;
\item there exists such a triangle;
\item $n$ has a prime divisor congruent to $1\pmod4$.
\end{enumerate}
\end{theorem}

\begin{proof}
(ii)$\Rightarrow$(i) is immediate.

(iii)$\Rightarrow$(ii). Let $q\equiv1\pmod4$ be a prime divisor of $n$, say
$n=qm$. By Fermat's theorem $q=a^2+b^2$ with positive integers $a\neq b$ ($q$ is odd).
Then $2m(a^2-b^2)$ (in absolute value), $4mab$ and $2mq=2n$ are the sides of a
nondegenerate right triangle with integer sides, since
$(a^2-b^2)^2+(2ab)^2=q^2$. Its circumradius is half its hypotenuse, namely $n$.

(i)$\Rightarrow$(iii), contrapositively. Suppose that $n$ has no prime divisor
congruent to $1\pmod4$. Then there is no non-diameter Gaussian side. Indeed, such a
side would give $s^2+u^2=(2n)^2$ with $s,u>0$; dividing by $\gcd(s,u)^2$ gives a
primitive Pythagorean triple whose hypotenuse $2n/\gcd(s,u)>1$ is odd and has all
its prime divisors congruent to $1\pmod 4$, and these divide $n$.
Hence Proposition~\ref{cor:nongaussian-polygons} shows that every admissible
polygon has an even number of sides, so no odd-sided polygon exists.
\end{proof}

\section{Semiprime radii}\label{sec:semiprime}
The regular hexagon gives $a(n)\ge6$. For semiprime $n$ we bound the number of
chords per unit angle within each squarefree class. Zero-sum relations among
component vectors (\S\ref{ssec:components}) give the bound $T/U$ for a block of
total $T$ in every ordinary class (\S\ref{ssec:ordinary}), where $U=\pi/6$. In
the Gaussian and Eisenstein classes (\S\S\ref{ssec:gauss}--\ref{ssec:eis}) the
extra units allow more chords per unit angle: we list all candidate half-angles
as integral combinations of one unit angle and the angles attached to the split
primes, give each candidate an integer weight $\Phi\ge1$ that is additive over
rational blocks, and obtain the densities $4/3$ and $3$. Explicit integer
polygons attain these densities. Class splitting adds the bounds over all
classes (\S\ref{ssec:mixtures}, where major arcs and polygons with a diameter
are also treated), and \S\ref{ssec:assembly} assembles the proof of
Theorem~\ref{thm:main}.

Set $H=\pi/2$ and $U=\pi/6$, so $H=3U$. Blocks, rational blocks and pure blocks
are as in Definition~\ref{def:block}; the total of a rational block is always positive.

\subsection{Component bounds and zero-sum atoms}\label{ssec:components}

\begin{lemma}[Uniform component bounds]\label{lem:components}
For semiprime $n$, the following bounds hold for each non-diameter chord in a fixed squarefree class:
\begin{enumerate}[nosep,label=(\roman*)]
\item $n=pq$, $p<q$ odd primes: $(v_p,v_q)\in\{-1,0,1\}^2$, with no nonzero component at $2$;
\item $n=p^2$, $p$ odd: $v_p\in\{-2,-1,0,1,2\}$, with no nonzero component at $2$;
\item $n=2p$, $p$ odd: $(v_2,v_p)\in\{-1,0,1\}^2$, with a coordinate that can be nonzero only when that prime splits;
\item $n=4$: only the prime $2$ can contribute, with $|v_2|\le2$;
\item $n=p$ prime (used only in Remark~\ref{rem:primes}): $v_p\in\{-1,0,1\}$ and no other component occurs.
\end{enumerate}
\end{lemma}
\begin{proof}
Let $s$ be a non-diameter side in class $d$, with $\kappa,A,B,M$ as in
\eqref{eq:primitive-norm}--\eqref{eq:norm-equation}, $K=\Q(\sqrt{-d})$ and
$z=A+B\sqrt{-d}$. Let $r$ be a rational prime that splits, $(r)=\mathfrak r\bar{\mathfrak r}$, and put
$e_1=\operatorname{ord}_{\mathfrak r}z$, $e_2=\operatorname{ord}_{\bar{\mathfrak r}}z$. By
\eqref{eq:idealnorm}, $e_1+e_2=2\operatorname{ord}_rM$ and $v_r=(e_1-e_2)/2$; in particular
$|v_r|\le\operatorname{ord}_rM$, and $v_r=0$ if $r\nmid M$.

\emph{Odd split primes.} Let $r$ be odd with $r\mid M$. If $e_1,e_2\ge1$, then
$z\in\mathfrak r\cap\bar{\mathfrak r}=r\OK$. Since $2\OK\subseteq\mathbb Z[\sqrt{-d}]$,
we get $2z\in r\,\mathbb Z[\sqrt{-d}]$, i.e.\ $r\mid2A$ and $r\mid2B$, hence $r\mid A,B$
because $r$ is odd; this contradicts $\gcd(A,B)=1$. Therefore $\min(e_1,e_2)=0$ and
$|v_r|=\operatorname{ord}_rM\le\operatorname{ord}_r(2n)$.

\emph{The prime $2$.} The prime $2$ splits in $K$ if and only if $d\equiv7\pmod8$.
Suppose it splits and $2\mid M$. Since $d$ is odd, $A^2+dB^2=M^2$ is even and
$\gcd(A,B)=1$, both $A,B$ are odd; then $A^2+dB^2\equiv1+7\equiv0\pmod8$, so
$4\mid M$. With $\eta=(1+\sqrt{-d})/2\in\OK$ we have
$z=(A-B)+2B\eta\in2\OK$, so $e_1,e_2\ge1$.

\emph{Conclusion.} Recall $M\mid2n$. In case (i), $\operatorname{ord}_pM,\operatorname{ord}_qM\le1$, and
$\operatorname{ord}_2M\le1$ excludes a component at $2$ (it would need $4\mid M$). In case
(ii), $\operatorname{ord}_pM\le2$ and again $2$ cannot contribute. In case (iii),
$|v_p|\le1$; if $2$ splits and $2\mid M$, then $4\mid M\mid4p$ forces
$\operatorname{ord}_2M=2$, so $e_1+e_2=4$ with $e_1,e_2\ge1$ and $|v_2|\le1$. In case (iv)
only $2$ divides $2n=8$; if it splits and $2\mid M$, then $\operatorname{ord}_2M\in\{2,3\}$,
$e_1+e_2\in\{4,6\}$ with $e_1,e_2\ge1$, and hence $|v_2|\le1$ resp.\ $\le2$.
In case (v), $|v_p|\le1$ for odd $p$ by the first argument and, as $\operatorname{ord}_2(2p)=1$, no component at $2$ occurs;
for $p=2$ we have $M\mid4$, and a component at $2$ forces $4\mid M$, $e_1+e_2=4$ and $|v_2|\le1$.
\end{proof}

\begin{lemma}[Zero-sum atom lengths]\label{lem:atoms}
A minimal nonempty zero-sum multiset of nonzero vectors from
$\{-1,0,1\}^2$ has at most four members. The bound is sharp.
A minimal zero sum of nonzero vectors in $\{-2,-1,1,2\}$ has at
most three; in $\{-1,1\}$ it has exactly two.
\end{lemma}
\begin{proof}
If an opposite pair occurs, minimality leaves precisely that pair.
Otherwise the support contains at most one of each of the four pairs
of opposite directions, so it has at most four directions. Suppose
four directions occur. After an integral isometry of the square we
may take the two diagonal directions to be $(1,1)$ and $(1,-1)$,
with respective positive multiplicities $c,d$. Their total is
$(c+d,c-d)$. The horizontal direction must then be $(-1,0)$
with multiplicity $c+d$; the vertical direction has multiplicity
$|c-d|$. Since all four directions occur, $c\ne d$.
If $c>d$, the multiset contains the zero-sum triple
$\{(1,1),(-1,0),(0,-1)\}$; if $d>c$, it contains
$\{(1,-1),(-1,0),(0,1)\}$. Either is a proper nonempty
zero-sum submultiset, a contradiction. Thus the support has at
most three directions.

With no opposite pair, a positive dependence cannot have just one
or two distinct directions. For three noncollinear directions
$u,v,w$, the unique primitive positive dependence has coefficients
the absolute determinants $|\det(v,w)|$, $|\det(w,u)|$,
$|\det(u,v)|$, divided by their gcd. Each determinant is at most
$2$. Its absolute value can equal $2$ only for a pair of diagonal
vectors. With no opposite pair, there can be at most two diagonal
directions in the support, so at most one coefficient can be $2$.
Consequently the coefficients sum to at most $2+1+1=4$.
The multiset $\{(1,0),(1,0),(-1,1),(-1,-1)\}$ shows sharpness.
In one dimension the minimal possibilities, up to a global
sign, are $\{1,-1\}$, $\{2,-2\}$, $\{2,-1,-1\}$ (and $\{-2,1,1\}$).
\end{proof}

\subsection{Ordinary classes: a universal density bound}\label{ssec:ordinary}
For $d\ne1,3$, the only units of $\mathbb Q(\sqrt{-d})$ are $\pm1$.
This converts the bound on zero-sum atom lengths into a uniform bound
on the number of chords per unit of angular measure.

\begin{lemma}[No ordinary class exceeds the hexagonal density]\label{lem:ordinary}
For any semiprime radius $n$ and any rational block of $k$ half-angles, pure in one squarefree class $d\ne1,3$ and with total $T\in\pi\Q$, one has
\begin{equation}\label{eq:ordinary} k\le T/U.\end{equation}
In particular, every such block has $T\ge H$.
\end{lemma}
\begin{proof}
By \eqref{eq:zero} the component vectors of the members sum to zero. Decompose
them into minimal zero-sum atoms (repeatedly split off a minimal nonempty zero-sum
submultiset), allowing zero vectors as singleton atoms. Let $w_1,\dots,w_m$ be the
numbers $w$ of the members of one atom. By Lemma~\ref{lem:component} their product is a root of unity of
$\Q(\sqrt{-d})$, hence $\pm1$, so $\sum_i\beta_i\in\pi\Z$, where $0<\beta_i<H$.
A singleton atom ($m=1$) would give $\beta_1\in\pi\Z$, which is impossible. For $m=2$ we
would get $\beta_1+\beta_2\in\pi\Z\cap(0,\pi)=\emptyset$. By Lemmas~\ref{lem:components}
and \ref{lem:atoms} every remaining atom has $m=3$ or $m=4$ members
(for $n=p^2$ or $n=4$ only $m=3$ occurs, and for prime $n$ no atom remains).
Its $\beta$-sum is a multiple of $\pi$ in $(0,mH)$, hence exactly $\pi$, and its
half-angle sum is therefore $mH-\pi=(m-2)H$. Since $m\le3(m-2)$ for $m=3,4$, this atom
satisfies $m\le(m-2)H/U$. Summing over the atoms gives \eqref{eq:ordinary}; each atom
contributes at least $H$ to $T$.
\end{proof}

\subsection{Gaussian classes and integral constructions}\label{ssec:gauss}

\paragraph{Gaussian data at a prime.}
For $d=1$, $\mathbb Z[i]$ is a PID with unit group $\mu_4$.
For each odd prime $r\equiv1\pmod4$, fix a Gaussian prime $\gamma_r$
of norm $r$, so $r=\gamma_r\overline{\gamma_r}$.
Unique factorization shows that every Gaussian integer of norm $r^2$
is a unit multiple of $r$, $\gamma_r^2$, or $\overline{\gamma_r}^2$.
The first orbit gives the four axial representations of $r^2$ as a sum of two squares, while the other
two give the eight non-axial ordered signed representations.
Since $\gamma_r^2$ is not axial, $2\gamma_r^2$ has even nonzero coordinates; taking absolute values and ordering them
gives integers $g_r,h_r>0$ with
\begin{equation}\label{eq:gaussnorm}
\begin{gathered}
 g_r^2+h_r^2=4r^2,\quad 0<g_r<h_r,\quad g_r,h_r\equiv0\pmod2,\\
 \xi_r=\arcsin\bigl(g_r/(2r)\bigr)\in(0,H/2).
\end{gathered}
\end{equation}
(Equality $g_r=h_r$ would force $g_r^2=2r^2$, which is impossible.) The two
complementary chord half-angles at radius $r$ are $\xi_r$ and $H-\xi_r$. Let
$w_{\xi_r}=e^{i(H-\xi_r)}=(g_r+ih_r)/(2r)$ be the number attached to the chord
of half-angle $\xi_r$. Here $\kappa=2$, $M=r$, and Lemma~\ref{lem:components} gives
$|v_r(w_{\xi_r})|=1$; fix the order of the two prime ideals over $r$ so that
$v_r(w_{\xi_r})=+1$.

\paragraph{Angle representation.}
Let $w$ be the number of any non-diameter Gaussian chord at radius $n$, and let $a_r$
be its components at the oriented split primes. The quotient
$w/\prod_r w_{\xi_r}^{a_r}$ has trivial fractional ideal because all
of its split-prime components vanish; hence it is a unit of
$\mathbb Z[i]$, i.e.\ a power of $i$. Comparing arguments gives
\begin{equation}\label{eq:gaussangles}
 \alpha=cH+\sum_r a_r \xi_r,\quad c\in\mathbb Z,\quad 0<\alpha<H,
\end{equation}
where $a_r$ are the split-prime components with these orientations. The vector $(a_r)$ is
the component vector of the chord, so it is additive. A chord with all $a_r=0$ would have $\beta$ equal to the argument
of a unit, which is impossible for $0<\beta<H$.

\paragraph{Two distinct split primes.}
Let $n=pq$ with distinct primes $p,q\equiv1\pmod4$ and write $x>y$ for the two angles
$\xi_p,\xi_q$ in decreasing order. (Since the radius is symmetric in its prime factors,
we label the primes by angle size rather than numerical size.) They cannot be equal, by
\eqref{eq:zero}. By Lemma~\ref{lem:components}, $(a,b)\in\{-1,0,1\}^2$. For each of the
eight nonzero pairs there is exactly one admissible $c$ in \eqref{eq:gaussangles},
and the complete candidate list, with the weight $\Phi=4c+2a+b$, is
\begin{center}
\small
\renewcommand{\arraystretch}{1.1}
\begin{tabular}{@{}cccc@{\qquad}cccc@{}}
\toprule
$(a,b)$&$c$&$\alpha$&$\Phi$&$(a,b)$&$c$&$\alpha$&$\Phi$\\
\midrule
$(1,0)$&$0$&$x$&$2$&$(1,1)$&$0$&$x+y$&$3$\\
$(-1,0)$&$1$&$H-x$&$2$&$(-1,-1)$&$1$&$H-x-y$&$1$\\
$(0,1)$&$0$&$y$&$1$&$(1,-1)$&$0$&$x-y$&$1$\\
$(0,-1)$&$1$&$H-y$&$3$&$(-1,1)$&$1$&$H-x+y$&$3$\\
\bottomrule
\end{tabular}
\end{center}
(Here $0<y<x<H/2$, so $0<x-y<H$ and $0<x+y<H$.) Thus $\Phi\ge1$ on every candidate
and $\Phi\le3$. For a rational block, \eqref{eq:zero} gives $\sum a=\sum b=0$, so the
total is $T=H\sum c$ and $\sum\Phi=4\sum c=4T/H=\tfrac43\,T/U$. Hence a block of
total $tH$ has at most $4t$ sides, and a $\pi$-block at most eight. Equality is
attained by doubling the $H$-block
\begin{equation}\label{eq:gaussconstruct}
 (x-y)+y+y+(H-x-y)=H.
\end{equation}
To verify actual integer chord lengths, label the primes so that
$x=\xi_p>y=\xi_q$ (this labeling need not have $p<q$). With $n=pq$, the three distinct
half-angles in \eqref{eq:gaussconstruct} correspond respectively to the chord lengths
\begin{equation}\label{eq:gausslengths}
 \frac{|g_ph_q-h_pg_q|}{2},\qquad p g_q,\qquad
 \frac{h_ph_q-g_pg_q}{2}.
\end{equation}
(For instance $2n\sin y=pq\,g_q/q=pg_q$ and $2n\sin(H-x-y)=2n\cos(x+y)=(h_ph_q-g_pg_q)/2$.) The
final numerator is positive because $x+y<H$. All three are integers by
\eqref{eq:gaussnorm} (the numerators are divisible by $4$), and their
multiplicities in the octagon are $2,4,2$. Thus the eight-side
upper bound is attained whenever $p\ne q$ both split in $\mathbb Z[i]$.

\paragraph{A repeated split prime.}
For $n=p^2$ with $p\equiv1\pmod4$, put $x=\xi_p$. By Lemma~\ref{lem:components},
$a\in\{-2,-1,1,2\}$, and \eqref{eq:gaussangles} gives exactly four candidates:
\begin{center}
\small
\renewcommand{\arraystretch}{1.1}
\begin{tabular}{@{}ccccc@{}}
\toprule
$a$&$1$&$-1$&$2$&$-2$\\
\midrule
$c$&$0$&$1$&$0$&$1$\\
$\alpha$&$x$&$H-x$&$2x$&$H-2x$\\
$\Phi=3c+a$&$1$&$2$&$2$&$1$\\
\bottomrule
\end{tabular}
\end{center}
(here $0<2x<H$). So $\Phi\in\{1,2\}$, $\sum\Phi=3\sum c=3T/H=T/U$, and a $\pi$-block has at most six sides.

\paragraph{Exactly one split prime.}
Suppose that exactly one odd prime $r$ splits in $\mathbb Z[i]$ and that it divides $n$ to
exponent one, as for $n=2r$ or $n=rq$ with $q\equiv3\pmod4$ (this includes the
one-split-prime case of odd semiprimes). The only nonzero component is
$a=\pm1$. Putting $x=\xi_r$, the complete candidate set is $\{x,H-x\}$, i.e.\
$(c,a)=(0,1)$ and $(1,-1)$, and both candidates have value $1$ under the weight
$\Phi=2c+a$. Hence a rational block with total $T$ has exactly $\sum\Phi=2T/H$ members,
that is, at most $\tfrac23\,T/U\le T/U$. If no odd prime splits in $\mathbb Z[i]$ there is no
non-diameter Gaussian chord at all.

\subsection{Eisenstein classes and integral constructions}\label{ssec:eis}
For two distinct split primes the list of candidate angles depends on whether
$x+y<U$ or $x+y>U$ (notation below); only the first case allows eighteen sides.

Let $K=\mathbb Q(\sqrt{-3})$, $\mathcal O_K^\times=\mu_6$, $H=3U$.
Its ring of integers is a PID. Every $\nu\in\OK$ has the form $\nu=(s+u\sqrt{-3})/2$ with
integers $s\equiv u\pmod 2$; conversely, any integer solution of
$s^2+3u^2=4r^2$ satisfies $s\equiv u\pmod2$ (reduce modulo $4$), so it defines such a $\nu$ with
$\Norm(\nu)=r^2$. A solution with $s,u>0$ corresponds to $0<\arg\nu<H$, and then
$\tan\arg\nu=u\sqrt3/s$ is the $\tan\beta$ of \eqref{eq:tanbeta}.

\begin{lemma}[Eisenstein chords at a prime radius]\label{lem:eisdata}
Let $r\equiv1\pmod3$ be prime. The positive solutions of $s^2+3u^2=4r^2$ are
$(s,u)=(r,r)$, with half-angle $U$, and exactly three further solutions, whose
half-angles are
\[
 x_r,\qquad 2U-x_r,\qquad 2U+x_r\qquad\text{for a unique }x_r\in(0,U).
\]
In particular there is exactly one solution with $0<s<r$, namely the one with half-angle
$x_r=\arcsin(s_r/(2r))$, and it has the smallest $s$ among all nontrivial
solutions. This proves the existence and uniqueness claimed before Theorem~\ref{thm:main}.
\end{lemma}
\begin{proof}
Because $\OK$ is a unique-factorization domain and
$r=\varpi\bar\varpi$ for a prime $\varpi$, every $\nu$ of norm $r^2$
is associated to exactly one of $r$, $\varpi^2$, or $\bar\varpi^2$.
Unit multiples of $r$ give only the half-angle $U$ in the open sector
(namely $\arg\nu=2U$, i.e.\ $(s,u)=(r,r)$).

For the two nontrivial orbits, let
$\phi\in(0,2U)$ be $2\arg\varpi$ reduced modulo $2U$.
Neither $\phi=0$ nor $\phi=U$ is possible: either would make
$\varpi^2/r=\varpi/\bar\varpi$ a root of unity, contradicting
its nontrivial fractional ideal $(\varpi)/(\bar\varpi)$.
Thus $\delta=\min\{\phi,2U-\phi\}$ belongs to $(0,U)$.
Applying the six units, conjugating, and retaining arguments in
$(0,H)=(0,3U)$ gives exactly three arguments,
$\beta=\delta,\ 2U-\delta,\ 2U+\delta$.
Their chord half-angles $\alpha=H-\beta$ are, in increasing order,
$U-\delta,\ U+\delta,\ 3U-\delta$.
Putting $x_r=U-\delta$ gives precisely
$\{x_r,2U-x_r,2U+x_r\}$, where $0<x_r<U$.
Thus exactly one nontrivial chord half-angle is below $U$; since
$s=2r\sin\alpha$, it is the unique solution with $0<s<r$ and it has the
smallest $s$ among the nontrivial solutions.
\end{proof}
The choice of $\varpi$ versus $\bar\varpi$ merely interchanges the
nontrivial orbits. For a prime $r\equiv1\pmod3$ we write $(s_r,u_r)$ for the solution of
Lemma~\ref{lem:eisdata} with half-angle $x_r$; thus
\begin{equation}\label{eq:eisnorm}
\begin{aligned}
 s_r^2+3u_r^2&=4r^2,\quad s_r,u_r>0,\quad 0<s_r<r,\\
 s_r&\equiv u_r\pmod2,\quad x_r=\arcsin(s_r/(2r)).
\end{aligned}
\end{equation}
Let $w_{x_r}=e^{i(H-x_r)}$ be the number of the chord of half-angle $x_r$ at radius $r$.
Then $r\nmid\gcd(s_r,u_r)$ (otherwise $(s_r/r,u_r/r)$ would be a positive solution of $a^2+3b^2=4$, i.e.\ the trivial solution $(1,1)$), so $r\mid M$ and
$|v_r(w_{x_r})|=\operatorname{ord}_rM=1$ by the proof of Lemma~\ref{lem:components}. We orient the prime-ideal pair
over $r$ so that $v_r(w_{x_r})=+1$.

\paragraph{Angle representation.}
Let $p\ne q$ be distinct primes $\equiv1\pmod3$ and write $x>y$ for $x_p,x_q$ in
decreasing order (their independence follows from \eqref{eq:zero}).
For an arbitrary Eisenstein chord $w$ at radius $pq$, write its oriented
split-prime component vector $(a,b)\in\{-1,0,1\}^2$ and let
$w_x=e^{i(H-x)}$, $w_y=e^{i(H-y)}$. By the component
identity \eqref{eq:zero}, the quotient $w/(w_x^a w_y^b)$
is a unit, hence a sixth root of unity. Therefore for some
$j\in\mathbb Z$ its complementary angle satisfies
$\beta\equiv a(H-x)+b(H-y)+2jU\pmod{2\pi}$.
Subtracting from $H=3U$ and absorbing multiples of $2\pi=12U$ into $j$ gives
$\alpha=cU+ax+by$ with $c=3-3a-3b-2j$;
so $c\equiv1-a-b\pmod 2$.
This proves that every possible non-diameter Eisenstein
half-angle lies in the candidate set
\begin{equation}\label{eq:eisenparity}
\begin{aligned}
 \alpha&=cU+ax+by,\quad (a,b)\in\{-1,0,1\}^2,\\
 &c\in\mathbb Z,\quad c\equiv1-a-b\pmod2,\quad 0<\alpha<3U.
\end{aligned}
\end{equation}
The candidate list will be shown to be exact by
Lemma~\ref{lem:special-realization}.

\paragraph{Complete two-prime Eisenstein catalogue.}
Here $0<y<x<U$, so $0<x-y<U$ and $0<x+y<2U$. The equality $x+y=U$ is impossible:
the half-angles $x_p$, $x_q$ and the half-angle $U$ of a hexagon side (whose
component vector is zero) would satisfy $x_p+x_q-U=0\in\pi\Q$ with component
vector $(1,1)\ne0$, contradicting \eqref{eq:zero}. Imposing the parity rule and $0<cU+ax+by<3U$ for each of the nine pairs $(a,b)$
gives the table below, which lists the admissible values of $c$ in each of the two chambers
$x+y<U$ and $x+y>U$. In square brackets is the weight of the candidate: $\Phi_+=3c+a+b$ in the
first chamber and $\Phi_-=c+a+b$ in the second.
\begin{center}
\small
\renewcommand{\arraystretch}{1.1}
\begin{tabular}{@{}ccc@{}}
\toprule
$(a,b)$&chamber $x+y<U$&chamber $x+y>U$\\
&entries $c\,[\Phi_+]$&entries $c\,[\Phi_-]$\\
\midrule
$(-1,-1)$&$1\,[1],\ 3\,[7]$&$3\,[1]$\\
$(-1,1)$&$1\,[3],\ 3\,[9]$&$1\,[1],\ 3\,[3]$\\
$(0,0)$&$1\,[3]$&$1\,[1]$\\
$(1,-1)$&$1\,[3]$&$1\,[1]$\\
$(1,1)$&$1\,[5]$&$-1\,[1],\ 1\,[3]$\\
$(-1,0)$&$2\,[5]$&$2\,[1]$\\
$(0,-1)$&$2\,[5]$&$2\,[1]$\\
$(0,1)$&$0\,[1],\ 2\,[7]$&$0\,[1],\ 2\,[3]$\\
$(1,0)$&$0\,[1],\ 2\,[7]$&$0\,[1],\ 2\,[3]$\\
\bottomrule
\end{tabular}
\end{center}
(Here $(0,0)$ is the hexagonal angle $U$, and $(a,b)=(0,1),(1,0)$ with $c=0$ are $y$ and $x$.)
Each chamber has thirteen candidates. The weights in the first chamber are
$1,1,1,3,3,3,5,5,5,7,7,7,9$ and those in the second are $1^{(9)},3^{(4)}$ (nine ones and four threes).
In particular $\Phi_+\ge1$ with maximum $9$, and $\Phi_-\ge1$ with maximum $3$. For a rational block,
$\sum a=\sum b=0$, so the total is $T=U\sum c$ and
\[
 \sum\Phi_+=3\sum c=3\,T/U,\qquad \sum\Phi_-=\sum c=T/U .
\]
Hence $k\le3T/U$ if $x+y<U$ and $k\le T/U$ if $x+y>U$.

\paragraph{The eighteen-gon.}
In the active case $x=x_p$, $y=x_q$ (in either order) and $x+y<U$, six copies of
\begin{equation}\label{eq:eisconstruct}
 x+y+(U-x-y)=U
\end{equation}
form an $18$-gon of radius $n=pq$. The first two chord lengths
are $q s_p$ and $p s_q$, and the third is
\begin{equation}\label{eq:eislengths}
 S=\frac{3u_pu_q-s_ps_q-3(s_pu_q+u_ps_q)}{4}.
\end{equation}
Indeed $2pq\sin x_p=qs_p$, $2pq\sin x_q=ps_q$, and, using
$\sin x_r=s_r/(2r)$, $\cos x_r=\sqrt3u_r/(2r)$,
\begin{multline*}
 \sin(U-x_p-x_q)=\tfrac12\cos(x_p+x_q)-\tfrac{\sqrt3}2\sin(x_p+x_q)\\
 =\frac{3u_pu_q-s_ps_q-3(s_pu_q+u_ps_q)}{8pq}.
\end{multline*}
So $S=2pq\sin(U-x_p-x_q)$, and it is positive because $U-x_p-x_q>0$. Set $R=4S$. Modulo $4$,
using $3\equiv-1\pmod4$, one has
\[
 R\equiv(s_p-u_p)(u_q-s_q)\equiv0\pmod4,
\]
since $s_r\equiv u_r\pmod2$. Thus $S\in\mathbb Z$,
and \eqref{eq:eisconstruct} gives exact closure with six
copies of each chord length.

\paragraph{One split prime, no split prime, and prime squares.}
If only one odd prime splits in $\mathbb Q(\sqrt{-3})$ (with exponent one in $n$), write $x=x_r$;
the same argument as above gives $a\in\{-1,0,1\}$, $c\equiv1-a\pmod2$ in
\eqref{eq:eisenparity}, and the complete candidate set is
$\{U,\ x,\ 2U-x,\ 2U+x\}$, i.e.\ $(c,a)=(1,0),(0,1),(2,-1),(2,1)$ with weights
$\Phi=c+a=1,1,1,3$. A rational block therefore has $k\le\sum\Phi=T/U$. If no odd prime splits, the
only candidate is $\alpha=U$ (the hexagonal angle), again with $k\le T/U$.

For a repeated split prime $n=p^2$ with $x=x_p$, the coefficient
$a$ lies in $\{-2,-1,0,1,2\}$. If $2x<U$, its complete
candidate list is
\begin{equation}\label{eq:eisrepeatedplus}
 (c,a)\in\{(1,-2),(3,-2),(2,-1),(1,0),
 (0,1),(2,1),(1,2)\};
\end{equation}
the weights $\Phi_+=3c+a$, in the displayed order, are
$(1,7,5,3,1,7,5)$, all at least $1$.
If $2x>U$, the candidates are
\begin{equation}\label{eq:eisrepeatedminus}
 (c,a)\in\{(3,-2),(2,-1),(1,0),(0,1),(2,1),(-1,2),(1,2)\};
\end{equation}
the weights $\Phi_-=c+a$, in the displayed order, are
$(1,1,1,1,3,1,3)$, all at least $1$.
Therefore the respective $\pi$-block bounds are $18$ and $6$ (because $\sum\Phi_+=3T/U$ and
$\sum\Phi_-=T/U$), with maximal weights $7$ and $3$.
The boundary $2x=U$ is impossible, since
$\sin^2x=s_p^2/(4p^2)$ is rational while
$\sin^2(\pi/12)=(2-\sqrt3)/4$ is not.
When $2x<U$, six copies of $x+x+(U-2x)=U$
form an $18$-gon: the first two equal chord lengths are $p s_p$,
while the third is
\begin{equation}\label{eq:eisrepeatedlength}
 S=\frac{2p^2-s_p^2-3s_pu_p}{2}.
\end{equation}
This is a positive integer because $s_p\equiv u_p\pmod2$;
thus the side multiset is $\{(p s_p)^{[12]},S^{[6]}\}$.

\begin{lemma}[Realizability of special-field catalogues]\label{lem:special-realization}
For a semiprime radius, every candidate half-angle in the Gaussian or
Eisenstein catalogue of this section satisfying its stated component,
unit-phase and open-sector conditions is realized by an integer chord.
Thus these two catalogues are exact, not merely upper-bound lists.
\end{lemma}
\begin{proof}
In either field $K=\Q(i)$ or $K=\Q(\sqrt{-3})$, each split odd prime
$r$ has the normalized prime-radius direction
$w_r=\eta_r/r$, where $\eta_r\in\mathcal O_K$ and
$|w_r|=1$. For an admissible component vector $(a_r)$ the catalogue fixes an
integer unit shift. In the Gaussian case this is the integer $c$
in $\alpha=cH+\sum_r a_r\xi_r$, corresponding to a power of $i$.
In the Eisenstein case the parity condition
$c\equiv1-\sum_r a_r\pmod2$ ensures that
$j=(3-3\sum_r a_r-c)/2$ is an integer, so the prescribed
phase is a sixth root of unity. This applies also to the
repeated-prime case, with $|a_r|\le\nu_r(n)=2$.
Thus the chosen unit $\varepsilon\in\mathcal O_K^\times$ exists,
and we may form
\[
 w=\varepsilon\prod_{r}w_r^{a_r}.
\]
For negative $a_r$ use $w_r^{-1}=\overline{w_r}$; hence
$w_r^{a_r}=\eta_{r,a_r}/r^{|a_r|}$ with
$\eta_{r,a_r}\in\mathcal O_K$.
Because $|a_r|\le\nu_r(n)$ for every candidate, multiplication by
$n$ clears all denominators, and $nw\in\mathcal O_K$.
The phase and sector restrictions yield $w=e^{i(H-\alpha)}$ with
$0<\alpha<H$, so both real and imaginary parts of $2nw$ are positive.
In the Gaussian case $2nw=s+iu$ with positive integers $s,u$.
In the Eisenstein case
$2nw=s+u\sqrt{-3}$ with positive integers $s,u$ of the same parity,
because $\mathcal O_K=\Z[(1+\sqrt{-3})/2]$.
In either case $s^2+du^2=(2n)^2$, giving $0<s<2n$ and a chord of
integer length $s$ with half-angle $\alpha$.
The coordinate-zero Eisenstein candidate $\alpha=U$ is provided by the
regular hexagon. The presence of inert or ramified rational factors of
$n$ does not affect this construction, since they are already included
in the multiplier $n$.
\end{proof}

\subsection{Mixtures, major arcs and diameters}\label{ssec:mixtures}
The catalogues are summarized in Table~\ref{tab:weights}. The two last columns are the only information about them
that is needed below: the \emph{density} $\rho=\sum\Phi\big/(T/U)$ controls
minor-arc polygons, and the maximal weight controls major-arc relations.

\begin{table}[htbp]
\centering
\small
\renewcommand{\arraystretch}{1.15}
\begin{tabular}{@{}>{\raggedright\arraybackslash}p{0.50\textwidth}ccc@{}}
\toprule
Class and case ($n$ semiprime) & weight $\Phi$ & density $\rho$ & $\max\Phi$\\
\midrule
Gaussian, two distinct split primes & $4c+2a+b$ & $4/3$ & $3$\\
Gaussian, $n=p^2$, $p\equiv1\ (4)$ & $3c+a$ & $1$ & $2$\\
Gaussian, one split prime & $2c+a$ & $2/3$ & $1$\\
Eisenstein, two distinct split primes, $x+y<U$ & $3c+a+b$ & $3$ & $9$\\
Eisenstein, two distinct split primes, $x+y>U$ & $c+a+b$ & $1$ & $3$\\
Eisenstein, $n=p^2$, $p\equiv1\ (3)$, $2x<U$ & $3c+a$ & $3$ & $7$\\
Eisenstein, $n=p^2$, $p\equiv1\ (3)$, $2x>U$ & $c+a$ & $1$ & $3$\\
Eisenstein, one split prime & $c+a$ & $1$ & $3$\\
Eisenstein, no split prime & $c$ & $1$ & $1$\\
Ordinary ($d\ne1,3$) & (Lemma~\ref{lem:ordinary}) & $1$ & ---\\
\bottomrule
\end{tabular}
\caption{Weights on the candidate half-angles. For a rational block with total $T$,
$k\le\sum\Phi=\rho\,T/U$, since every candidate has $\Phi\ge1$ and $\sum a=\sum b=0$.
``Two distinct split primes'' means $n=pq$ with $p\ne q$ both split in the given field; ``one split
prime'' means that exactly one prime divisor of $n$ splits, to exponent one.}
\label{tab:weights}
\end{table}

\begin{proposition}[Uniform semiprime density bounds]\label{prop:semiprime-density}
Let $n$ be semiprime. Every rational block of $k$ half-angles, pure in one squarefree class $d$ and with total $T$, satisfies
$k\le\rho_d T/U$, where $U=\pi/6$ and $\rho_d$ is the
following upper-bound coefficient (not necessarily the exact density
in Table~\ref{tab:weights}):
\begin{equation}\label{eq:densities}
\rho_d=\begin{cases}
3,& d=3,\ \text{two distinct split odd primes with }x+y<U,\\
  &\phantom{d=3,}\text{or one repeated split odd prime with }2x<U,\\
4/3,&d=1,\ \text{two distinct split odd primes},\\
1,&\text{otherwise}.
\end{cases}
\end{equation}
Consequently, any minor-arc polygon without diameter at radius $n$ has at most
$6\max_d\rho_d$ sides.
\end{proposition}
\begin{proof}
For ordinary classes the bound is Lemma~\ref{lem:ordinary}. For the Gaussian and Eisenstein classes it is the
column $\rho$ of Table~\ref{tab:weights}, which was derived in \S\ref{ssec:gauss} and \S\ref{ssec:eis} (the
one-split-prime Gaussian density $2/3$ is bounded by $1$). Let $P$ be a minor-arc polygon without diameter; its
sides form a $\pi$-block. By Lemma~\ref{lem:splitting}, the restriction of this $\pi$-block to each
nonempty squarefree class is a rational block, with total $T_d>0$ and $k_d$ members, say; then
$k_d\le\rho_dT_d/U$ and $\sum_dT_d=\pi$. Summing over the classes gives
$k\le(\max_d\rho_d)\pi/U=6\max_d\rho_d$.
\end{proof}

\begin{lemma}[No mixed non-diameter major arcs]\label{lem:nomix}
In a non-diameter major-arc equation
$\alpha_*=\sum_{j=1}^{k-1}\alpha_j<H$, all sources and the target
belong to the same squarefree class.
\end{lemma}
\begin{proof}
Let $d_*$ be the target class, and suppose the foreign sources (those not in $d_*$) have
total angle $E$ with $0<E<H$ (each source is positive and $\sum\alpha_j=\alpha_*<H$). Apply
Lemma~\ref{lem:unitlattice} to $\alpha_*-\sum_j\alpha_j=0$ with signs $+1$ for the target and $-1$ for the sources.
If $d_*\ne3$, the signed subtotal in $d_*$ lies in $H\Z$ (parts (i), (ii) of the lemma); but it equals
$\alpha_*-\sum_{j\in d_*}\alpha_j=E\in(0,H)$, a contradiction. If $d_*=3$, every foreign class $e\ne3$ containing source angles
has signed subtotal $-\sum_{j:\,d_j=e}\alpha_j\in H\Z$.
Its \emph{positive} source-angle sum therefore belongs to $H\Z$
and is at least $H$. Thus $E\ge H$, again impossible.
\end{proof}

\begin{lemma}[Uniform major-arc bounds]\label{lem:major}
A non-diameter major-arc polygon at a semiprime radius has at most
four sides in an ordinary or Gaussian class, at most ten in the
active two-prime Eisenstein case, at most eight in the active
repeated-prime Eisenstein case, and at most four in any other
Eisenstein case.
\end{lemma}
\begin{proof}
By Lemma~\ref{lem:nomix} all sides lie in one class. Let $k$ be the number of sides, so there are $k-1\ge2$ sources.

\emph{Ordinary class.} The source and target component vectors
are nonzero: a zero component would make the complementary angle
the argument of a unit, impossible for $0<\beta<H$.
In the augmented signed zero sum
$\{-v_*,v_1,\ldots,v_{k-1}\}$, no proper nonempty source-only
submultiset can have sum zero, since Lemma~\ref{lem:ordinary}
would give it an angle sum at least $H$, whereas
$\sum\alpha_j=\alpha_*<H$. A proper zero-sum submultiset
containing $-v_*$ would leave such a source-only submultiset
as its nonempty complement. Hence the augmented multiset
is a minimal zero sum, and Lemma~\ref{lem:atoms}
gives $k\le4$ (and $k\le3$ for repeated-prime
or $n=4$ cases).

\emph{Gaussian and Eisenstein classes.} The component equation
is $a_*=\sum a_j$, $b_*=\sum b_j$ (where present), by \eqref{eq:zero}.
The genuine major-arc equation is
$c_*U+a_*x+b_*y=\sum_j(c_jU+a_jx+b_jy)$,
so, after subtracting the component identities, $(c_*-\sum_jc_j)U=0$ and hence $c_*=\sum_j c_j$ (with $H$ replacing $U$ in the Gaussian
case).
Thus each weight $\Phi$ of Table~\ref{tab:weights} is \emph{exactly additive}: $\Phi(\alpha_*)=\sum_j\Phi(\alpha_j)$.
Every source contributes at least $1$, so $k-1\le\Phi(\alpha_*)\le\max\Phi$. By
Table~\ref{tab:weights}, $\max\Phi=3$ in the two-prime Gaussian case (so $k\le4$), $2$ in the
repeated-prime Gaussian case ($k\le3$), and $1$ in the one-prime Gaussian case (so
there is no major-arc polygon with at least two sources). The Eisenstein maxima are $9$ (two-prime active),
$7$ (repeated active), and $3$ in all other cases; hence $k\le10,8,4$, respectively.
\end{proof}

\begin{lemma}[Polygons with a diameter]\label{lem:diameter}
Let $n$ be semiprime and $\rho_{\max}=\max_d\rho_d$ over the classes present. A polygon at radius $n$ with a diameter has at most
$1+3\rho_{\max}$ sides, that is, at most $10$, $5$, or $4$ sides when
$\rho_{\max}=3$, $4/3$, or $1$, respectively. These numbers are strictly below the minor-arc
bounds $18$, $8$, $6$.
\end{lemma}
\begin{proof}
The other sides form an $H$-block (Remark~\ref{rem:diameter}). Exactly as in the proof of
Proposition~\ref{prop:semiprime-density}, with total $H=3U$ in place of $\pi=6U$, it has at most $3\rho_{\max}$
members, i.e.\ at most $9$, $4$, or $3$.
\end{proof}

\subsection{Assembly: proof of Theorem~\ref{thm:main}}\label{ssec:assembly}
\paragraph{Small radii $2p$ and $4$.}
For $n=2p$, $p$ odd, all classes satisfy the baseline density
$\rho\le1$: ordinary classes by Lemma~\ref{lem:ordinary}, and the Gaussian and Eisenstein classes by the one-split-prime
and no-split-prime rows of Table~\ref{tab:weights} (the prime $2$ neither splits in $\mathbb Q(i)$ nor in $\mathbb Q(\sqrt{-3})$).
Hence minor-arc polygons have at most $6$ sides; by Lemmas~\ref{lem:major} and \ref{lem:diameter}, non-diameter major-arc
polygons and polygons with a diameter have at most four sides. Therefore $a(2p)\le6$, and the regular hexagon attains six.
For $n=4$, the squarefree classes of the non-diameter
lengths $s=1,\ldots,7$ are, in order,
$7,15,55,3,39,7,15$ (from $4n^2-s^2=64-s^2=63,60,55,48,39,28,15$). None is Gaussian, and the only
Eisenstein half-angle is $U$ (at $s=4$). Every ordinary
class has only the single possible split component at $2$,
with $|v_2|\le2$, so Lemma~\ref{lem:ordinary} applies and gives density $1$;
ordinary major-arc polygons have at most three sides (Lemma~\ref{lem:major}),
while no major-arc relation can arise from
angles all equal to $U$ (it would read $U=rU$ with $r\ge2$). Thus again $a(4)=6$.

\paragraph{Odd semiprimes.}
Let $n=pq$ with $p\le q$ odd primes. By Proposition~\ref{prop:semiprime-density}
and Table~\ref{tab:weights}, the largest class density $\rho_{\max}$ is
\begin{itemize}[nosep,leftmargin=2.6em]
\item $3$, if $p\equiv q\equiv1\pmod3$ and $x_p+x_q<U$ (for $p=q$ this reads $2x_p<U$);
\item $4/3$, if this fails but $p<q$ and $p\equiv q\equiv1\pmod4$;
\item $1$, otherwise.
\end{itemize}
Indeed, the Eisenstein class has density $3$ only in the two ``active''
cases of Table~\ref{tab:weights}; the Gaussian density exceeds $1$ only for two
distinct split primes; and the Eisenstein case takes precedence when both
conditions hold, because its density is larger. The equalities $x_p+x_q=U$ and
$2x_p=U$ cannot occur (see \S\ref{ssec:eis}), so the first alternative fails
exactly when $x_p+x_q>U$, or when one of $p,q$ is not $\equiv1\pmod3$.

\begin{proof}[Proof of Theorem~\ref{thm:main}]
\emph{Upper bounds.} By Proposition~\ref{prop:semiprime-density}, a minor-arc
polygon without diameter has at most $6\rho_{\max}$ sides, that is, at most
$18$, $8$ or $6$ sides according to the three alternatives above. For
$n=2p$ and $n=4$ the bound is $6$ by the first paragraph of this subsection. By
Lemmas~\ref{lem:geometry}, \ref{lem:nomix}, \ref{lem:major}
and~\ref{lem:diameter}, neither major-arc polygons nor polygons with a diameter
exceed these maxima.

\emph{Lower bounds.} In the active Eisenstein cases the integral
constructions \eqref{eq:eisconstruct} and
\eqref{eq:eisrepeatedlength} give $18$ sides; otherwise
two distinct Gaussian split primes give the integral
$8$-side construction \eqref{eq:gaussconstruct}. The regular
hexagon yields six in all remaining cases, including $n=4$ and $n=2p$.

\emph{The integer criterion.} By the computation after \eqref{eq:eislengths},
\[
\sin(U-x_p-x_q)=
\frac{3u_pu_q-s_ps_q-3(s_pu_q+u_ps_q)}{8pq}.
\]
Since $0<x_p+x_q<2U$, this sine is positive precisely when
$x_p+x_q<U$; this is equivalent to the claimed integer
criterion \eqref{eq:integer-criterion}. The repeated-prime
case is obtained by putting $p=q$.

\emph{The special values.} If $p=2$ or $p=3$, then $p\not\equiv1\pmod3$ and $p\not\equiv1\pmod4$, so neither exceptional
case occurs; hence $a(2q)=a(3q)=6$ for every prime $q$, since the exceptional cases require \emph{both} primes to be $\equiv1\pmod3$ or both $\equiv1\pmod4$ (and $a(pq)$ is symmetric in $p,q$).
\end{proof}

\begin{corollary}[Rigidity of exceptional extremisers]\label{cor:unique-extremisers}
For a semiprime radius $n=pq$, if $a(n)=18$ or $a(n)=8$, the
side-length multiset of an extremal polygon is unique up to ordering.
For $p\ne q$ in the active Eisenstein case, it is
$\{(q s_p)^{[6]},(p s_q)^{[6]},S^{[6]}\}$ with $S$ given by
\eqref{eq:eislengths}; in the active repeated-prime case it is
$\{(p s_p)^{[12]},S^{[6]}\}$ with $S$ given by
\eqref{eq:eisrepeatedlength}. In the Gaussian eight-side case the three
lengths in \eqref{eq:gausslengths} occur with multiplicities $2,4,2$.
\end{corollary}
\begin{proof}
The bounds in Lemmas~\ref{lem:major} and~\ref{lem:diameter} are strictly
below the exceptional maxima, so an extremiser is a minor-arc polygon
without a diameter. Indeed, for every class with $\rho_d<\rho_{\max}$ the contribution is
$k_d\le\rho_d T_d/U<\rho_{\max}T_d/U$ whenever $T_d>0$.
Equality in the global bound therefore forces all sides to be in the
unique class of maximal density: Eisenstein for $18$ and Gaussian for $8$.
Every angle in that class must have weight $\Phi=1$, since
$k\le\sum\Phi$ and $k=\sum\Phi$ at equality.

If $p\ne q$ and $a(pq)=18$, the only Eisenstein weight-one candidates in
the active chamber are $x$, $y$ and $U-x-y$, with component vectors
$(1,0)$, $(0,1)$ and $(-1,-1)$. If their multiplicities are $A,B,C$,
component cancellation gives $A=B=C$, and their total is $CU=\pi$,
so $A=B=C=6$. If $p=q$, the weight-one angles are $x$ and $U-2x$,
with component values $1$ and $-2$; hence $A=2C$ and
$Ax+C(U-2x)=CU=\pi$, yielding $A=12$, $C=6$.

In the Gaussian case $a(pq)=8$, the weight-one angles are
$x-y$, $y$, and $H-x-y$, with components $(1,-1)$, $(0,1)$,
and $(-1,-1)$. The zero-sum equations give multiplicities
$A=C$ and $B=2A$. Their half-angle total is $AH=\pi$,
so $(A,B,C)=(2,4,2)$. By the component criterion,
these listed angles cannot coincide: an equality of two distinct
candidates would give a rational angle relation with a nonzero
component vector. The component vector associated with each chord is determined uniquely
by its prime-ideal valuations; there is no ambiguity from alternative
catalogue encodings. Finally, $\alpha_n(s)$ is strictly increasing in
$s$, so uniqueness of angles implies uniqueness of chord lengths.
Different cyclic orders are not asserted to give congruent polygons.
\end{proof}

\begin{remark}[Primes and $n=1$]\label{rem:primes}
For $n=1$ the only chords are $s=1$ ($\alpha=U$) and the diameter $s=2$, so $a(1)=6$. For a prime $n=p$,
Lemma~\ref{lem:components}(v) gives $1$-dimensional component vectors with entries in $\{-1,0,1\}$, whose minimal zero sums have two members
only (Lemma~\ref{lem:atoms}); so no ordinary class contains a rational block, and (Gaussian one-split-prime and
Eisenstein rows of Table~\ref{tab:weights}) every class has density at most $1$. The arguments above then give $a(p)=6$
(minor-arc polygons have at most $6$ sides; a diameter plus an $H$-block has at most $4$; major arcs at most $4$).
\end{remark}
\section{Growth along fixed powers}\label{sec:powers}

Increasing the exponents at fixed prime support behaves differently from the
semiprime case: the following construction gives unbounded side counts along
$N^e$ for every fixed $N\ge2$. Write
\begin{equation}\label{eq:Ddv}
 4N^2-1=dv^2,
\end{equation}
where $d$ is squarefree and $v\ge1$. Put
\[
 z=1+v\sqrt{-d},
 \qquad \Norm(z)=4N^2,
\]
and write
\[
 z^j=x_j+y_j\sqrt{-d}\qquad(j\ge1).
\]

\begin{lemma}[The quadratic chord tower]\label{lem:chord-tower}
For every $j\ge1$, both $x_j$ and $y_j$ are divisible by $2^{j-1}$. If
\[
 X_j=\frac{x_j}{2^{j-1}},\qquad Y_j=\frac{y_j}{2^{j-1}},
\]
then
\begin{equation}\label{eq:tower-norm}
 X_j^2+dY_j^2=4N^{2j}.
\end{equation}
In particular, $|X_j|$ is an integer chord length at radius $N^j$.
Moreover, if
\[
 \theta=\arg z\in(0,\pi/2),
 \qquad \tau=\frac{2\theta}{\pi},
\]
then $\tau$ is irrational and the associated half-angle is
\begin{equation}\label{eq:folded-alpha}
 \lambda_j:=\alpha_{N^j}(|X_j|)
=\arcsin|\cos(j\theta)|.
\end{equation}
\end{lemma}

\begin{proof}
The real parts satisfy the trace recursion
\[
 x_{j+2}=2x_{j+1}-4N^2x_j,
 \qquad x_0=x_1=1.
\]
With $X_j=x_j/2^{j-1}$ this becomes
\begin{equation}\label{eq:Xrec}
 X_{j+2}=X_{j+1}-N^2X_j,
 \qquad X_1=1,\quad X_2=1-2N^2,
\end{equation}
so all $X_j$ are integers. The imaginary parts satisfy the same trace recursion, with
$y_1=v$ and $y_2=2v$; after division by $2^{j-1}$ the same integral
recurrence applies. Hence $2^{j-1}\mid y_j$ (with $y_0=0$). Dividing
$\Norm(z^j)=(4N^2)^j$ by $4^{j-1}$ proves
\eqref{eq:tower-norm}.

Since $z/(2N)=e^{i\theta}$,
\[
 \frac{X_j}{2N^j}=\cos(j\theta),
\]
which gives \eqref{eq:folded-alpha}. If $\theta/\pi$ were rational,
then $e^{i\theta}$ would be a root of unity and
\[
 2\cos\theta=\frac1N
\]
would be a rational algebraic integer. It would therefore be an integer,
which is impossible for $N\ge2$. Hence $\theta/\pi$ and $\tau$ are
irrational. Consequently $X_j\ne0$, so this is a positive integer chord
length.
\end{proof}

For $e\ge j$, scaling the chord in Lemma~\ref{lem:chord-tower} gives the
integer side
\begin{equation}\label{eq:scaled-chord}
 s_{e,j}=N^{e-j}|X_j|
\end{equation}
at radius $N^e$, with the same half-angle $\lambda_j$.

We use one elementary fact about the Stern--Brocot path of an irrational
number in $(0,1)$.

\begin{lemma}[Odd--odd bracketing edges]\label{lem:odd-odd}
Let $\tau\in(0,1)$ be irrational. Along its Stern--Brocot path there are
infinitely many bracketing edges
\begin{equation}\label{eq:odd-edge}
 \frac{a_1}{b_1}<\tau<\frac{a_2}{b_2},
 \qquad b_1a_2-a_1b_2=1,
\end{equation}
 for which both numerators $a_1$ and $a_2$ are odd. The denominators of these
edges are unbounded.
\end{lemma}

\begin{proof}
For neighbors $a_1/b_1<a_2/b_2$ with $b_1a_2-a_1b_2=1$ the numerators $a_1,a_2$
are not both even. Call such an edge \emph{mixed} if exactly one numerator is
odd. The mediant $(a_1+a_2)/(b_1+b_2)$ has odd numerator if the edge is mixed
and even numerator if it is odd--odd. The Stern--Brocot path of $\tau$ starts
at the mixed edge $0/1<\tau<1/1$. A step from a mixed edge replaces either the
even endpoint, which gives an odd--odd edge, or the odd endpoint, which gives
again a mixed edge with the same even endpoint $\varepsilon$. If the second
alternative occurred at every step from some point on, the other endpoint would
converge to $\varepsilon$ and $\tau=\varepsilon$ would be rational. Hence every
mixed edge is followed later by an odd--odd edge, and every odd--odd edge is
followed immediately by a mixed edge. So odd--odd edges occur infinitely
often. The edges are pairwise distinct and only finitely many fractions have
bounded denominator, so the denominators are unbounded.
\end{proof}

\begin{lemma}[Balanced odd--odd brackets]\label{lem:balanced-brackets}
Enumerate the odd--odd bracketing edges of Lemma~\ref{lem:odd-odd}
in their order along the Stern--Brocot path, and let
\[
 r_j=\frac{\min(b_{1,j},b_{2,j})}{\max(b_{1,j},b_{2,j})}.
\]
Then $\limsup_{j\to\infty}r_j\ge\sqrt2-1$.
Moreover $e_j=\max(b_{1,j},b_{2,j})$ increases strictly.
\end{lemma}
\begin{proof}
Let an odd--odd edge have denominators $x,y$ and put $E=x+y$.
Its mediant has even numerator and denominator $E$, so the next
edge is mixed. Its even endpoint has denominator $E$, and its odd
endpoint has denominator $R\in\{x,y\}$, corresponding to the
endpoint of the old edge that was retained. At a mixed edge the
mediant has odd numerator. Hence the next step either replaces the
odd endpoint, leaving the even endpoint fixed and producing another
mixed edge, or replaces the even endpoint, producing an odd--odd
edge. Let $k\ge0$ be the number of steps of the first kind before
the first step of the second kind; such a step eventually occurs by
Lemma~\ref{lem:odd-odd}. After those $k$ steps the odd endpoint has
denominator $R+kE$. Replacing the even endpoint by the mediant
therefore gives the next odd--odd edge, with denominators
\[
 R+kE,\qquad R+(k+1)E.
\]
The larger denominator strictly exceeds $\max(x,y)$. For $k\ge1$,
\[
 \frac{R+kE}{R+(k+1)E}>\frac12
 \quad\Longleftrightarrow\quad R+(k-1)E>0,
\]
so the ratio of the smaller denominator to the larger one exceeds
$1/2$. If $k=0$ and $r=\min(x,y)/\max(x,y)$, the new ratio is
$A(r)=r/(1+2r)$ when $R=\min(x,y)$, and
$B(r)=1/(2+r)$ when $R=\max(x,y)$.

If $k\ge1$ occurs infinitely often, the conclusion follows because
$1/2>\sqrt2-1$. Otherwise, eventually $k=0$. Then $B$ must occur
infinitely often. Indeed, if all sufficiently late transitions
were of type $A$, the smaller-denominator odd endpoint would be
retained at every transition: an $A$ step with $k=0$ produces
new denominators $R<R+E$, so the retained endpoint remains the
smaller one at the next transition. Hence one fixed rational
fraction $a/R$ would be an endpoint of every sufficiently late
bracketing edge. If $a'/b'$ is the other endpoint, neighboring
fractions satisfy $|a/R-a'/b'|=1/(Rb')$, and consequently
\[
 |\tau-a/R|<\frac1{Rb'}\longrightarrow0,
\]
since the other denominator $b'$ strictly increases. This forces
$\tau=a/R\in\mathbb Q$, contradicting the irrationality of $\tau$.
Set $c=\sqrt2-1$, so $B(c)=c$. If $\limsup r_j<c$, then
some $\varepsilon>0$ satisfies $r_j\le c-\varepsilon$ eventually.
At every later occurrence of $B$ we would have
\[
 r_{j+1}=\frac1{2+r_j}\ge\frac1{2+c-\varepsilon}>c,
\]
a contradiction.
\end{proof}

\begin{proof}[Proof of Theorem~\ref{thm:power-growth}]
Choose an odd--odd edge \eqref{eq:odd-edge}. Since $b_1a_2-a_1b_2=1$,
\[
 0<b_1\tau-a_1<\frac1{b_2}\le1,
 \qquad
 0<a_2-b_2\tau<\frac1{b_1}\le1.
\]
Since $a_1$ and $a_2$ are odd and both distances are less than $1$ (odd integers
are $2$ apart), these are the distances of $b_1\tau$ and $b_2\tau$ from the
nearest odd integers $a_1$ and $a_2$. (If $m$ is odd and $|t-m|\le1$, then
$|\cos(Ht)|=|\sin(H(t-m))|$, and $\arcsin|\cos(Ht)|=H|t-m|$.) Formula
\eqref{eq:folded-alpha} therefore gives, with $H=\pi/2$,
\[
 \lambda_{b_1}=H(b_1\tau-a_1),
 \qquad
 \lambda_{b_2}=H(a_2-b_2\tau).
\]
The determinant identity now yields
\begin{equation}\label{eq:half-block-growth}
 b_2\lambda_{b_1}+b_1\lambda_{b_2}
 =H(b_1a_2-a_1b_2)=H.
\end{equation}
Doubling \eqref{eq:half-block-growth} gives a $\pi$-block with
$2(b_1+b_2)$ sides. By \eqref{eq:scaled-chord}, both required chord lengths are
integers at the common radius $N^e$, where $e=\max\{b_1,b_2\}$. Hence
\begin{equation}\label{eq:growth-lower}
 a(N^e)\ge2(b_1+b_2)=2e(1+r),
 \qquad r=\frac{\min(b_1,b_2)}{\max(b_1,b_2)}.
\end{equation}
By Lemmas~\ref{lem:odd-odd} and~\ref{lem:balanced-brackets},
the exponents $e_j$ increase without bound, and
$\limsup_j r_j\ge\sqrt2-1$. Consequently
$\limsup_{e\to\infty}a(N^e)/e\ge2\sqrt2$.
The same construction shows that $a(N^e)$ is unbounded. By Proposition~\ref{prop:scaling}, $a(N^e)$ is
nondecreasing in $e$. Hence $a(N^e)\to\infty$.
\end{proof}

The following example makes the construction explicit; the proof above does
not depend on it.

\begin{example}[A thirty-eight-gon at radius $8192$]\label{ex:8192}
Take $N=2$. Then $d=15$, $v=1$, and $z=1+\sqrt{-15}$. For
\[
 \tau=\frac{2\arg z}{\pi}
\]
one has the odd--odd Farey edge
\[
 \frac56<\tau<\frac{11}{13},
 \qquad 6\cdot11-5\cdot13=1.
\]
The Farey inequalities can be certified without trigonometric
approximation. Writing $\theta=\arg z=\arccos(1/4)$, we have
$5/6<\tau$ because
$\cos(5\pi/12)=(\sqrt6-\sqrt2)/4>1/4$;
the final inequality follows by squaring from $48<49$.
Likewise $\tau<11/13$ since
$\cos(11\pi/26)=\sin(\pi/13)<\pi/13<1/4$.
For $N=2$, the integer recurrence \eqref{eq:Xrec} begins
\[
 1,-7,-11,17,61,-7,-251,-223,781,1673,-1451,-8143,-2339,
\]
so in particular the scaled power coefficients are exactly
\[
 \frac{\Re(z^6)}{2^5}=X_6=-7,
 \qquad
 \frac{\Re(z^{13})}{2^{12}}=X_{13}=-2339.
\]
At the common radius $2^{13}=8192$, the first chord scales by $2^7$ and
has length $896$, while the second has length $2339$. Thus
\[
 26\alpha_{8192}(896)+12\alpha_{8192}(2339)=\pi,
\]
so
\[
 \{896^{[26]},2339^{[12]}\}
\]
is a thirty-eight-gon and $a(8192)\ge38$.
\end{example}

\begin{figure}[htbp]
\centering
\begin{tikzpicture}[x=1cm,y=1cm]
\draw[->,thick] (-0.3,0)--(7.5,0);
\foreach \v/\x in {0.830/0,0.835/1.75,0.840/3.5,0.845/5.25,0.850/7.0}{
  \draw (\x,0.07)--(\x,-0.07); \node[font=\scriptsize,below=2pt] at (\x,-0.07) {$\v$};}
\def\xa{1.1667}\def\xt{3.1986}\def\xb{5.6538}
\draw[line width=2.2pt] (\xa,0)--(\xb,0);
\foreach \x in {\xa,\xt,\xb}{\fill (\x,0) circle (2.2pt);}
\node[font=\small,above=5pt] at (\xa,0) {$\frac56$};
\node[font=\small,above=5pt] at (\xt,0) {$\tau$};
\node[font=\small,above=5pt] at (\xb,0) {$\frac{11}{13}$};
\draw[decorate,decoration={brace,amplitude=3pt,mirror}] (\xa,-0.62)--(\xt,-0.62);
\draw[decorate,decoration={brace,amplitude=3pt,mirror}] (\xt,-0.62)--(\xb,-0.62);
\node[font=\scriptsize,anchor=north east] at (\xt-0.05,-0.95) {$\tau-\tfrac56=\dfrac{\lambda_6}{6H}$};
\node[font=\scriptsize,anchor=north west] at (\xt+0.05,-0.95) {$\tfrac{11}{13}-\tau=\dfrac{\lambda_{13}}{13H}$};
\end{tikzpicture}
\caption{The Farey bracket $\frac56<\tau<\frac{11}{13}$ of Example~\ref{ex:8192},
drawn to scale; $\tau=2\arg(1+\sqrt{-15})/\pi\approx0.83914$. Here $H=\pi/2$ and
$\lambda_j=\alpha_{2^j}(|X_j|)$. The distances from $\tau$ to the endpoints are
$\lambda_6/(6H)$ and $\lambda_{13}/(13H)$, and $6\cdot11-5\cdot13=1$ gives
$13\lambda_6+6\lambda_{13}=H$.
The plotted positions were rendered with AI-assisted TikZ code from these values.}
\label{fig:farey-example}
\end{figure}
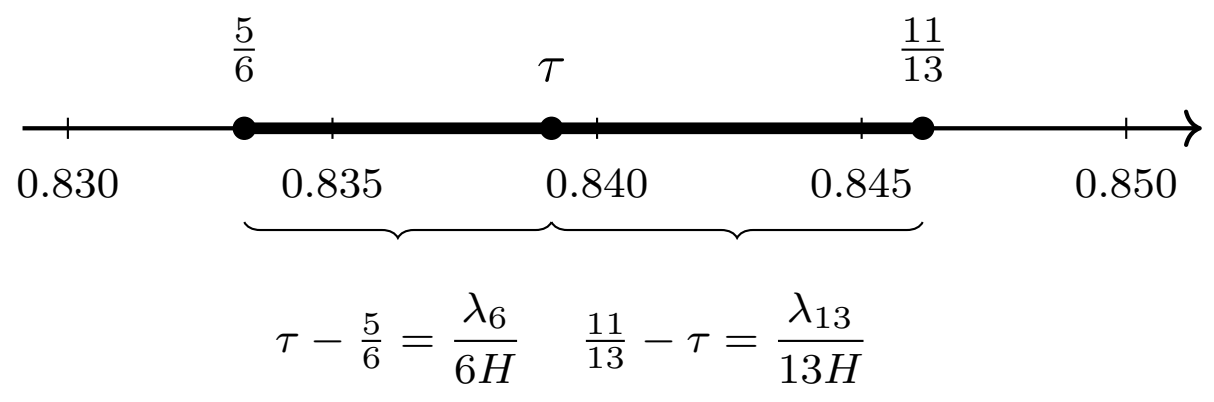

\subsection{Polynomial upper bounds along powers}\label{ssec:poly-upper}
The component bound also yields uniform control in fixed prime support.
Write $N=\prod_{p\mid N}p^{c_p}$, let $r=\omega(N)$, and set
\begin{equation}\label{eq:atom-capacity}
 L_e=\prod_{p\mid N}(2rec_p+1).
\end{equation}

\begin{lemma}[Components at arbitrary radii]\label{lem:components-all}
For every non-diameter chord at radius $n$ in any squarefree class,
$|v_p(w)|\le\nu_p(n)$ for every split prime $p$.
All components at primes not dividing $n$ vanish.
\end{lemma}
\begin{proof}
Use the notation $A,B,M,z$ from \eqref{eq:primitive-norm}--
\eqref{eq:norm-equation}. Let $p$ be an odd split prime, and write $(p)=\mathfrak p\bar{\mathfrak p}$.
If both prime ideals divided $(z)$, then
$z\in\mathfrak p\cap\bar{\mathfrak p}=p\mathcal O_K$.
Since $2\mathcal O_K\subseteq\Z[\sqrt{-d}]$, this would give
$p\mid2A$ and $p\mid2B$, hence $p\mid A,B$, contradicting
$\gcd(A,B)=1$. Thus at most one of the two prime-ideal valuations
of $z$ is positive, and
$|v_p(w)|=\nu_p(M)\le\nu_p(n)$ because $M\mid2n$
and $p$ is odd.

For $p=2$, splitting requires $d\equiv7\pmod8$. If $2\mid M$,
then the equation $A^2+dB^2=M^2$ and $\gcd(A,B)=1$ force
$A,B$ to be odd. Consequently $A^2+dB^2\equiv1+7\equiv0\pmod8$,
so $4\mid M$. With
$\eta=(1+\sqrt{-d})/2\in\mathcal O_K$ we have
$z=(A-B)+2B\eta\in2\mathcal O_K$.
If $t=\nu_2(M)$, the two valuations at the primes over $2$ are
at least $1$ and sum to $2t$, so
$|v_2(w)|\le t-1\le\nu_2(n)$ since $M\mid2n$.
If $2\nmid M$ the component vanishes. In particular, the
coordinate at $2$ vanishes when $2\nmid n$.
\end{proof}

\begin{lemma}[Size of a zero-sum atom]\label{lem:atom-capacity}
For $n=N^e$, every minimal nonempty zero-sum multiset of
component vectors from a fixed class has at most $L_e$ members.
\end{lemma}
\begin{proof}
Pad vectors by zeros to one coordinate for every prime divisor of
$N$, including the primes inert or ramified in the fixed class.
By Lemma~\ref{lem:components-all}, all coordinates are integers and
$|v_p|\le ec_p$. A zero vector is a singleton atom.
Otherwise apply the Steinitz rearrangement lemma
\cite[Lemma~1]{KleinSteinitz} to the zero-sum sequence with the
scaled norm $\|(x_p)\|=\max_{p\mid N}|x_p|/(ec_p)$.
There is an ordering of its $m$ terms with all partial sums in the
box $|x_p|\le rec_p$. The $m$ integer partial sums
$S_0,S_1,\ldots,S_{m-1}$ are distinct: a repetition would give a
proper nonempty zero-sum submultiset, contradicting minimality.
The box contains exactly $L_e$ integer points, so $m\le L_e$.
\end{proof}

\begin{theorem}[Polynomial bounds along powers]\label{thm:polynomial-bound}
For any fixed integer $N\ge2$ and every $e\ge1$,
\begin{equation}\label{eq:polynomial-bound}
 a(N^e)\le6\prod_{p\mid N}(2\omega(N)e\nu_p(N)+1).
\end{equation}
In particular $a(N^e)=O_N(e^{\omega(N)})$.
For every prime $p$, the stronger estimate
\begin{equation}\label{eq:linear-prime-upper}
 a(p^e)\le12e
\end{equation}
holds. If $q\equiv1\pmod4$ is prime, then
\begin{equation}\label{eq:theta-prime}
 2e<a(q^e)\le12e,\qquad a(q^e)=\Theta(e).
\end{equation}
\end{theorem}
\begin{proof}
Consider first a minor-arc polygon without diameter at radius
$n=N^e$. The sum of its half-angles is $\pi=6U$.
By Lemma~\ref{lem:splitting} its angles split into rational blocks,
one for each squarefree class. Decompose the component zero sums in
each block into minimal nonempty atoms. For an atom, the component sum is zero; by the reverse implication of
Lemma~\ref{lem:component} its angle sum belongs to $\pi\Q$,
and Lemma~\ref{lem:unitlattice} places this sum in $U\Z$.
Each atom contains only positive angles, so its total is at least $U$.
A zero-vector atom is a singleton; the only possible zero-component
non-diameter chord is the Eisenstein hexagon chord of angle $U$.
There are at most six atoms in all, and each has at most $L_e$
members by Lemma~\ref{lem:atom-capacity}; thus there are at most
$6L_e$ sides.

For a polygon with a diameter, its remaining angles total
$H=3U$. The same argument gives at most three atoms, hence
at most $1+3L_e\le6L_e$ sides.

For a major-arc polygon, Lemma~\ref{lem:nomix} applies at every
radius (its proof does not use the semiprime restriction).
All sides therefore belong to a single class. Decompose the
signed zero-sum multiset $\{-v_*,v_1,\ldots,v_{k-1}\}$ into
minimal atoms. Exactly one contains the negative target; all other
atoms consist entirely of source angles. These source-only atoms
have positive totals in $U\Z$ and their sum is at most
$\alpha_*<H=3U$, so there are at most two. Therefore
$k\le3L_e$. Together, the three cases prove
\eqref{eq:polynomial-bound}.

For $N=p$ there is only one potentially nonzero coordinate,
lying in $[-e,e]\cap\Z$. A minimal nonempty zero sum of such
integers has at most $2e$ members: for an atom without a zero
entry, start with a positive term and thereafter choose a term
opposite in sign to the current nonzero partial sum. This is
always possible because the remaining terms sum to its negative.
More precisely, the $m$ distinct partial sums
$S_0=0,S_1,\ldots,S_{m-1}$ all lie in $[1-e,e]\cap\Z$:
starting with a positive entry keeps every nonzero partial sum
at most $e$, and an opposite-sign entry can lower it no further
than $1-e$; $S_0$ belongs to this same interval. A repetition would
contradict minimality. As the interval contains exactly $2e$ integers,
the atom has at most $2e$ members. Thus one can replace $L_e$ by
$2e$ in the three polygon cases,
proving $a(p^e)\le12e$.

Finally, if $q\equiv1\pmod4$, Lemma~\ref{lem:gaussian-farey}
gives an integral Gaussian $\pi$-block with $2(b_1+b_2)$
sides at radius $q^e$. For Farey neighbors of order $e$ one has
$b_1+b_2>e$. Combining this lower bound with
\eqref{eq:linear-prime-upper} gives \eqref{eq:theta-prime}.
\end{proof}

\section{Parity with one Gaussian split prime}\label{sec:parity}
We show that $a(n)$ is even if at most one prime divisor of $n$ is congruent to
$1\pmod4$. An odd extremal polygon is first reduced to a pure Gaussian block;
the maximal length of such blocks is then computed by Farey neighbors.

Put $H=\pi/2$. A \emph{Gaussian $H$-block} is an $H$-block (Definition~\ref{def:block}) pure in the Gaussian class $d=1$,
and a \emph{Gaussian $\pi$-block} is defined similarly. Let $H_G(n)$ and $P_G(n)$
be their maximum cardinalities at radius $n$, with value $0$ if no such block exists
(for instance if the class is absent).

\begin{lemma}[Reduction of an odd maximum]\label{lem:odd-maximum-reduction}
If $a(n)$ is odd, then an extremal polygon is a pure Gaussian minor-arc
polygon without diameter. In particular, $a(n)\le P_G(n)$.
\end{lemma}

\begin{proof}
Let $P$ be an extremal polygon with $k=a(n)$ sides, $k$ odd. By
Proposition~\ref{cor:nongaussian-polygons}, $P$ must contain a non-diameter side of the
Gaussian class. We treat the three kinds of polygons of Remark~\ref{rem:diameter}.

\emph{Polygon with a diameter.} The remaining $k-1$ sides form an $H$-block with, say,
$r=k-1\ge2$ members. Doubling that block gives a $\pi$-block, hence by
Lemma~\ref{lem:geometry} a polygon with $2r>r+1=k$ sides, contradicting extremality.

\emph{Minor-arc polygon without diameter.} The sides form a $\pi$-block. Apply
Lemma~\ref{lem:unitlattice} (with all $\varepsilon_i=1$): the Gaussian subtotal $\mathcal A_1$ lies
in $H\Z$. If other squarefree classes occur, the Gaussian subtotal and the non-Gaussian
subtotal $\pi-\mathcal A_1\in H\Z$ are both positive, hence both equal $H$. Let their
side counts be $k_G$ and $k_N$, so $k=k_G+k_N$. Doubling either $H$-block gives a
$\pi$-block, hence an even polygon with $2k_G$ or $2k_N$ sides. Since $k$ is odd,
$k_G\ne k_N$, and
\[
 \max\{2k_G,2k_N\}>k_G+k_N=k,
\]
contradicting extremality. So no other class occurs: $P$ is purely Gaussian.

\emph{Major-arc polygon.} Here $\alpha_*=\sum_{j}\alpha_j<H$. Apply
Lemma~\ref{lem:unitlattice} to $\alpha_*-\sum_j\alpha_j=0$: the Gaussian subtotal lies in $H\Z$.
If the target is not Gaussian, this subtotal is $-\sum_{j\ \mathrm{Gaussian}}\alpha_j\in(-H,0]$,
hence zero; thus there is no Gaussian side, contradicting the first paragraph of this proof.
If the target is Gaussian, the Gaussian subtotal equals $\alpha_*-\sum_{j\ \mathrm{Gaussian}}\alpha_j=\sum_{j\ \mathrm{non\text{-}Gaussian}}\alpha_j\in[0,H)$, hence
zero, so every source side is Gaussian and the relation is purely Gaussian.
Write $4n^2-s_*^2=u_*^2$, where $s_*$ is the side with half-angle $\alpha_*$; the
complementary chord $u_*$ is an integer, is not a diameter, and has half-angle $H-\alpha_*$.
The source sides together with this complementary
chord form an $H$-block with $(k-1)+1=k$ members; doubling it gives a polygon with $2k>k$ sides,
a contradiction.

Hence $P$ is a pure Gaussian $\pi$-block without diameter, and so $k\le P_G(n)$.
\end{proof}

Assume now that $n$ has exactly one distinct prime divisor
$q\equiv1\pmod4$, and write $q^e\Vert n$. Choose a Gaussian prime
$\varpi$ with $\varpi\overline\varpi=q$ whose argument satisfies
\begin{equation}\label{eq:gaussian-farey-parameter}
 \vartheta=\arg\varpi\in(0,\pi/4);
 \qquad \sigma=\frac{4\vartheta}{\pi}\in(0,1)
\end{equation}
(of the eight numbers $u\varpi$, $u\overline\varpi$, $u\in\mu_4$, exactly one has argument in $(0,\pi/4)$).
The number $\sigma$ is irrational. Indeed, a rational value would make
$\varpi/\overline\varpi=e^{2i\vartheta}$ a root of unity in $\Q(i)$, i.e.\ a power of $i$; since
$0<2\vartheta<\pi/2$, this is impossible.

For $1\le j\le e$, let
\[
 \ell_j=\lfloor j\sigma\rfloor,
 \qquad f_j=j\sigma-\ell_j\in(0,1).
\]

\begin{lemma}[The rank-one Gaussian angle set]\label{lem:gaussian-angle-set}
The non-diameter Gaussian half-angles at radius $n$ are exactly
\begin{equation}\label{eq:gaussian-angle-list}
 Hf_j,
 \qquad H(1-f_j),
 \qquad 1\le j\le e.
\end{equation}
\end{lemma}

\begin{proof}
Let $s+iu$ be a Gaussian integer with $s^2+u^2=(2n)^2$. A prime
$r\equiv3\pmod4$ is inert in $\Z[i]$, so its exponent in $s+iu$ is
forced and contributes no argument. The prime $2$ is ramified and contributes
only a multiple of $H$, which can be absorbed into the unit. For the unique
split prime $q=\varpi\overline\varpi$, write the exponents of $\varpi$ and
$\overline\varpi$ in $s+iu$ as $a$ and $2e-a$. Their contribution to the
argument is
\[
 (a-(2e-a))\vartheta=2(a-e)\vartheta.
\]
Thus, modulo $H$, the argument of $s+iu$ is $\pm2j\vartheta=\pm Hj\sigma$ with
$j=|a-e|\in\{0,\dots,e\}$; the case $j=0$ would give an axial vector, that is a
diameter or a degenerate side, which is excluded. The half-angle of the chord is
$H$ minus this argument, and folding into $(0,H)$ gives
$Hf_j$ or $H(1-f_j)$ with $1\le j\le e$. Conversely, for $1\le j\le e$ the
Gaussian integer $\varpi^{e+j}\overline\varpi^{\,e-j}$, multiplied by the forced
factors at $2$ and at the primes $\equiv3\pmod4$ and by a unit, has norm
$(2n)^2$ and positive real and imaginary parts (it is not axial, because $\sigma$
is irrational), so it defines a chord; its half-angle is $H(1-f_j)$, and
exchanging $\varpi$ and $\overline\varpi$ gives $Hf_j$. None of these angles is
$0$ or $H$, since $\sigma$ is irrational.
\end{proof}

The primes $2$ and those congruent to $3\pmod4$ therefore change only the
rational scale of the corresponding chords; they create no new Gaussian
direction.

Let
\begin{equation}\label{eq:farey-rank-one}
 \frac{a_1}{b_1}<\sigma<\frac{a_2}{b_2}
\end{equation}
be the two neighbors of $\sigma$ in the Farey sequence of order $e$. Then
\[
 b_1a_2-a_1b_2=1,
 \qquad b_1,b_2\le e,
 \qquad b_1+b_2>e.
\]
Indeed, if $b_1+b_2\le e$, the mediant
$(a_1+a_2)/(b_1+b_2)$ would be a fraction of the Farey
sequence of order $e$ strictly between two consecutive terms.
Put
\[
 D=b_1+b_2,
 \qquad C=a_1+a_2.
\]

\begin{lemma}[Exact Gaussian block maxima]\label{lem:gaussian-farey}
With the notation above,
\begin{equation}\label{eq:PGHG}
 H_G(n)=D,
 \qquad P_G(n)=2D.
\end{equation}
\end{lemma}

\begin{proof}
Represent the two angle types in \eqref{eq:gaussian-angle-list} by
\[
 v_j=(-\ell_j,j),
 \qquad v'_j=(\ell_j+1,-j),
\]
so that a vector $(\mu,\nu)\in\Z^2$ represents the angle $H(\mu+\nu\sigma)$. Because $1$ and
$\sigma$ are linearly independent over $\Q$, an $H$-block
sums to $(1,0)$ and a $\pi$-block to $(2,0)$.

The mediant $C/D$ lies strictly between the Farey neighbors
\eqref{eq:farey-rank-one}. For every $1\le j\le e$, the function
$x\mapsto\lfloor jx\rfloor$ is constant on the \emph{open} interval
$(a_1/b_1,a_2/b_2)$: a jump in its interior would occur at a rational
with denominator at most $j\le e$, contradicting the adjacency of the
endpoints in the Farey sequence of order $e$. This assertion need not hold
at both endpoints of the closed Farey interval. The two values to which
we apply it, $\sigma$ and $C/D$, are strictly inside the open interval.
Moreover, $\gcd(C,D)=1$ by the determinant-one identity and $D>e\ge j$,
so $Cj/D$ is not an integer.
Consequently
\[
 \ell_j=\left\lfloor\frac{Cj}{D}\right\rfloor,
\]
and hence
\begin{equation}\label{eq:dual-bound}
 1\le Cj-D\ell_j\le D-1.
\end{equation}
Define the integral linear functional
\[
 F(\mu,\nu)=D\mu+C\nu.
\]
By \eqref{eq:dual-bound}, every available vector $v_j$ and $v'_j$ has
$1\le F\le D-1$ (indeed $F(v_j)=Cj-D\ell_j$ and $F(v'_j)=D-F(v_j)$). Since $F(1,0)=D$ and $F(2,0)=2D$, an $H$-block has at most $D$ sides and a
$\pi$-block at most $2D$ sides.

For the lower bound, Farey adjacency gives
\[
 b_2(-a_1,b_1)+b_1(a_2,-b_2)=(1,0).
\]
Here $(-a_1,b_1)=v_{b_1}$ and $(a_2,-b_2)=v'_{b_2}$, because
$a_1/b_1<\sigma<a_2/b_2$ gives $\ell_{b_1}=a_1$ and $\ell_{b_2}=a_2-1$. Hence $b_2$ copies of the first angle and $b_1$
copies of the second form an $H$-block with $b_1+b_2=D$ sides. Doubling it gives a
$\pi$-block with $2D$ sides. This proves \eqref{eq:PGHG}.
\end{proof}

\begin{proof}[Proof of Theorem~\ref{thm:parity}]
If $\omega_1(n)=0$, Theorem~\ref{thm:odd-polygon} (or
Proposition~\ref{cor:nongaussian-polygons}) shows that every
admissible polygon already has even side count. Suppose
$\omega_1(n)=1$. If $a(n)$ were odd, Lemma~\ref{lem:odd-maximum-reduction}
would give $a(n)\le P_G(n)$. But
Lemma~\ref{lem:gaussian-farey} shows that $P_G(n)=2D\ge4$ is attained by a Gaussian $\pi$-block, which is itself an admissible polygon, so
$a(n)\ge P_G(n)$. Thus an odd maximum would have to equal the even number
$P_G(n)$, a contradiction. Therefore $a(n)$ is even.
\end{proof}

\begin{corollary}[A further evenness criterion]\label{cor:gaussian-part-parity}
Let the Gaussian part of $n$ be
\[
 n_G=\prod_{\substack{q\mid n\\q\equiv1\ (\mathrm{mod}\ 4)}}q^{\nu_q(n)}.
\]
Then $P_G(n)=P_G(n_G)$. If $n_G$ is $1$, a prime power,
or the product of two distinct primes, then $a(n)$ is even.
\end{corollary}
\begin{proof}
In $\Z[i]$, each prime $q\equiv1\pmod4$ splits as
$q=\varpi_q\overline{\varpi_q}$. Writing
$e_q=\nu_q(n)$, the two exponents in a Gaussian integer
of norm $(2n)^2$ are $a_q$ and $2e_q-a_q$ with
$0\le a_q\le2e_q$. These are precisely the exponents
that determine its direction, modulo Gaussian units.
A prime congruent to $3$ modulo $4$ is inert and contributes only
a positive real factor; the ramified prime $2$ contributes only
a positive factor and a Gaussian unit (because $(1+i)^2=2i$).
Conversely, every admissible choice of the exponents $a_q$ and unit
phase occurs for both radii $n$ and $n_G$ after inserting the forced
factors. The resulting non-axial half-angle sets coincide, and hence
$P_G(n)=P_G(n_G)$.

The cases $n_G=1$ and $n_G=q^e$ follow from
Theorem~\ref{thm:parity}. If $n_G=pq$ with $p\ne q$ and
$p\equiv q\equiv1\pmod4$, apply the \emph{semiprime}
Gaussian two-prime catalogue in Section~\ref{ssec:gauss} at radius
$pq=n_G$, not at the possibly larger radius $n$.
It gives $P_G(pq)=8$: its density bound
gives at most eight sides, and \eqref{eq:gaussconstruct}
attains eight. If $a(n)$ were odd,
Lemma~\ref{lem:odd-maximum-reduction} would give
$a(n)\le P_G(n)=8$, while the Gaussian construction implies
$a(n)\ge P_G(n)=8$. This would make $a(n)=8$, a contradiction.
\end{proof}

\begin{remark}
The hypothesis $\omega_1(n)\le1$ is used in Lemma~\ref{lem:gaussian-angle-set},
where the Gaussian angles are indexed by a single exponent. For two
distinct split primes of exponent one, the Gaussian catalogue already
settles evenness by Corollary~\ref{cor:gaussian-part-parity}.
For general Gaussian parts containing multiple split primes, the
angles are indexed by a multidimensional box. We do not know
whether $a(n)$ is even for every $n$.
\end{remark}

\section{Examples and consequences}\label{sec:examples}

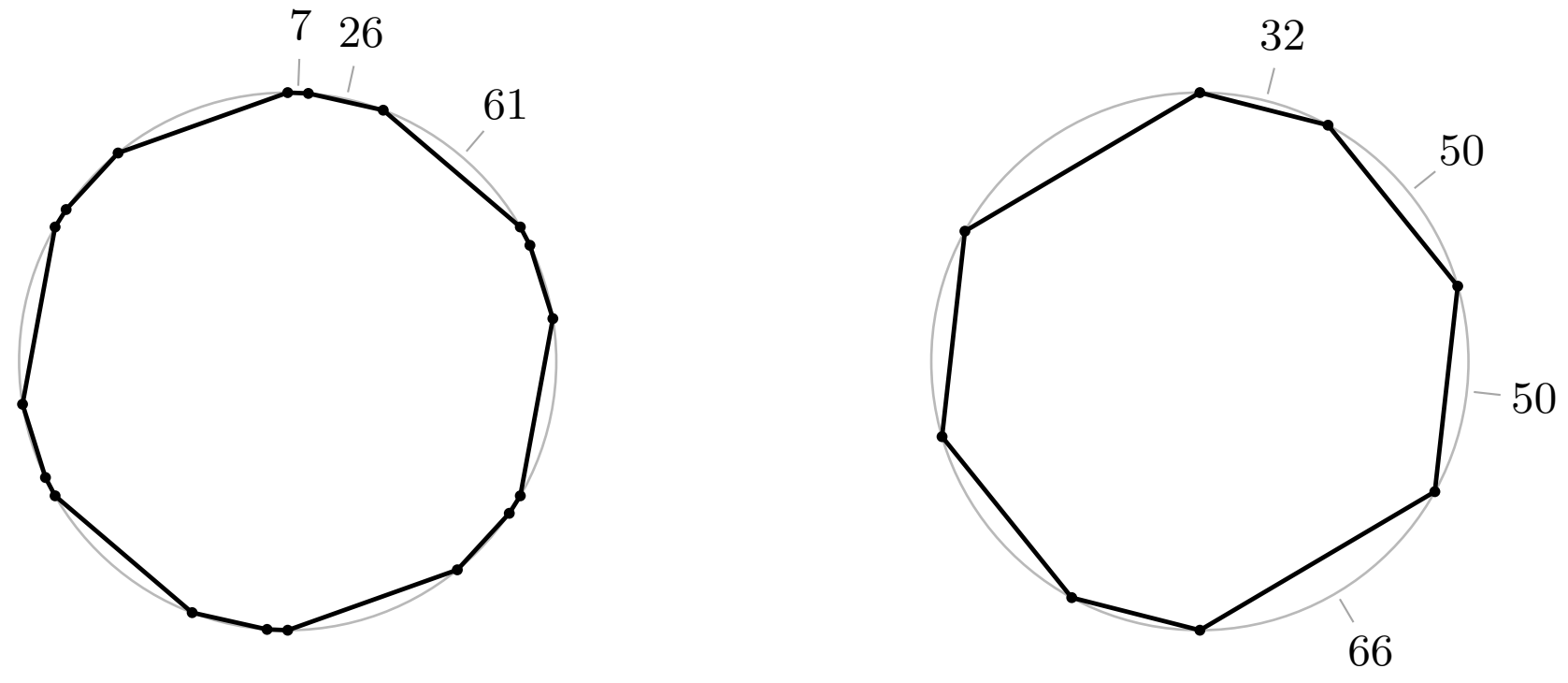
\begin{figure}[htbp]
\centering
\begin{minipage}[t]{0.49\textwidth}
\centering
\begin{tikzpicture}
\coordinate (V0) at (90.00000:1.90cm);
\coordinate (V1) at (85.59154:1.90cm);
\coordinate (V2) at (69.16512:1.90cm);
\coordinate (V3) at (30.00000:1.90cm);
\coordinate (V4) at (25.59154:1.90cm);
\coordinate (V5) at (9.16512:1.90cm);
\coordinate (V6) at (-30.00000:1.90cm);
\coordinate (V7) at (-34.40846:1.90cm);
\coordinate (V8) at (-50.83488:1.90cm);
\coordinate (V9) at (-90.00000:1.90cm);
\coordinate (V10) at (-94.40846:1.90cm);
\coordinate (V11) at (-110.83488:1.90cm);
\coordinate (V12) at (-150.00000:1.90cm);
\coordinate (V13) at (-154.40846:1.90cm);
\coordinate (V14) at (-170.83488:1.90cm);
\coordinate (V15) at (-210.00000:1.90cm);
\coordinate (V16) at (-214.40846:1.90cm);
\coordinate (V17) at (-230.83488:1.90cm);
\draw[gray!55,line width=0.5pt] (0,0) circle (1.90cm);
\draw[line width=0.9pt] (V0)--(V1)--(V2)--(V3)--(V4)--(V5)--(V6)--(V7)--(V8)--(V9)--(V10)--(V11)--(V12)--(V13)--(V14)--(V15)--(V16)--(V17)--cycle;
\foreach \k in {0,...,17} \fill (V\k) circle (1.1pt);
\draw[gray!70,line width=0.4pt] (87.7958:1.95cm)--(87.7958:2.14cm);
\node[font=\footnotesize,inner sep=1pt] at (87.7958:2.38cm) {$7$};
\draw[gray!70,line width=0.4pt] (77.3783:1.95cm)--(77.3783:2.14cm);
\node[font=\footnotesize,inner sep=1pt] at (77.3783:2.38cm) {$26$};
\draw[gray!70,line width=0.4pt] (49.5826:1.95cm)--(49.5826:2.14cm);
\node[font=\footnotesize,inner sep=1pt] at (49.5826:2.38cm) {$61$};
\path (-2.65,-2.6) rectangle (2.65,2.65);
\end{tikzpicture}

\smallskip
{\small (a) Radius $91$: sides $7^{[6]},26^{[6]},61^{[6]}$.}
\end{minipage}\hfill
\begin{minipage}[t]{0.49\textwidth}
\centering
\begin{tikzpicture}
\coordinate (V0) at (90.00000:1.90cm);
\coordinate (V1) at (61.49993:1.90cm);
\coordinate (V2) at (16.26020:1.90cm);
\coordinate (V3) at (-28.97953:1.90cm);
\coordinate (V4) at (-90.00000:1.90cm);
\coordinate (V5) at (-118.50007:1.90cm);
\coordinate (V6) at (-163.73980:1.90cm);
\coordinate (V7) at (-208.97953:1.90cm);
\draw[gray!55,line width=0.5pt] (0,0) circle (1.90cm);
\draw[line width=0.9pt] (V0)--(V1)--(V2)--(V3)--(V4)--(V5)--(V6)--(V7)--cycle;
\foreach \k in {0,...,7} \fill (V\k) circle (1.1pt);
\draw[gray!70,line width=0.4pt] (75.7500:1.95cm)--(75.7500:2.14cm);
\node[font=\footnotesize,inner sep=1pt] at (75.7500:2.38cm) {$32$};
\draw[gray!70,line width=0.4pt] (38.8801:1.95cm)--(38.8801:2.14cm);
\node[font=\footnotesize,inner sep=1pt] at (38.8801:2.38cm) {$50$};
\draw[gray!70,line width=0.4pt] (-6.3597:1.95cm)--(-6.3597:2.14cm);
\node[font=\footnotesize,inner sep=1pt] at (-6.3597:2.38cm) {$50$};
\draw[gray!70,line width=0.4pt] (-59.4898:1.95cm)--(-59.4898:2.14cm);
\node[font=\footnotesize,inner sep=1pt] at (-59.4898:2.38cm) {$66$};
\path (-2.65,-2.6) rectangle (2.65,2.65);
\end{tikzpicture}

\smallskip
{\small (b) Radius $65$: sides $32^{[2]},50^{[4]},66^{[2]}$.}
\end{minipage}
\caption{The polygons of Examples~\ref{ex:eis} and~\ref{ex:gauss}, not drawn to
the same scale. Vertices are placed according to the central angles
$2\alpha_n(s)$ of the consecutive sides, which repeat the pattern $7,26,61$
(six times) in (a) and $32,50,50,66$ (twice) in (b); one period is labelled.
The diagrams were rendered with AI-assisted TikZ code from these chord lengths.}
\label{fig:semiprime-examples}
\end{figure}

Figure~\ref{fig:semiprime-examples} shows the polygons of the following two
examples, which realize the lower bounds of Theorem~\ref{thm:main}.

\begin{example}[Eisenstein jump]\label{ex:eis}
For $p=7$ and $q=13$,
\[
 (s_7,u_7)=(2,8),
 \qquad (s_{13},u_{13})=(1,15).
\]
The first two chord lengths are $qs_p=26$ and $ps_q=7$, and formula \eqref{eq:eislengths} gives $S=61$. Thus
\[
 \{7^{[6]},26^{[6]},61^{[6]}\}
\]
is the side multiset of an eighteen-gon of radius $91$; hence $a(91)=18$ by
Theorem~\ref{thm:main}.
\end{example}

\begin{example}[Gaussian jump]\label{ex:gauss}
For $p=5$ and $q=13$, take $(g_5,h_5)=(6,8)$ and
$(g_{13},h_{13})=(10,24)$ from \eqref{eq:gaussnorm}; here $x=\xi_5>y=\xi_{13}$.
Formula~\eqref{eq:gausslengths} gives the integer side lengths
$|6\cdot24-8\cdot10|/2=32$, $5\cdot10=50$, and
$(8\cdot24-6\cdot10)/2=66$. Doubling the $H$-block
\eqref{eq:gaussconstruct}
therefore gives the side multiset
\[
 \{32^{[2]},50^{[4]},66^{[2]}\},
\]
so $a(65)=8$ by Theorem~\ref{thm:main}.
\end{example}

\FloatBarrier

\Needspace{10\baselineskip}
\begin{corollary}\label{cor:prime-squares}
For a prime $p$,
\[
 a(p^2)=18
\]
exactly when $p\equiv1\pmod3$ and $2x_p<\pi/6$; otherwise $a(p^2)=6$.
In particular,
\[
 a(49)=18,\qquad a(169)=18,\qquad a(361)=6.
\]
\end{corollary}
\begin{proof}
The prime-square criterion follows from Theorem~\ref{thm:main}. For the three stated examples, the positive solutions of
$s^2+3u^2=4p^2$ at $p=7,13,19$ are, respectively,
\[
\begin{aligned}
 p=7:&\quad (2,8),(7,7),(11,5),(13,3),\\
 p=13:&\quad (1,15),(13,13),(22,8),(23,7),\\
 p=19:&\quad (11,21),(19,19),(26,16),(37,5).
\end{aligned}
\]
Their completeness follows from Lemma~\ref{lem:eisdata} (one trivial solution $(p,p)$ and three
nontrivial ones). Hence the small-angle
representatives are $(s_p,u_p)=(2,8),(1,15),(11,21)$.
Since $3<(7/4)^2$ and $3>(12/7)^2$,
\[
 \frac1{16}<\sin^2(\pi/12)=\frac{2-\sqrt3}{4}<\frac1{14},
\]
the rational inequalities
$1/49,1/676<1/16$ and $121/1444>1/14$
show exactly that $2x_7,2x_{13}<\pi/6<2x_{19}$.
\end{proof}
By Theorems~\ref{thm:power-growth} and~\ref{thm:parity}, for every prime $p$ the
numbers $a(p^e)$ are even and tend to infinity with $e$.

\section{Open questions}\label{sec:outlook}

By Theorem~\ref{thm:parity} and Corollary~\ref{cor:gaussian-part-parity},
$a(n)$ is even when its Gaussian part is a prime power or the product
of two distinct primes (and also when it is $1$). By
Theorem~\ref{thm:main}, it is even for every semiprime $n$.
We do not know
whether
\[
 a(n)\equiv0\pmod2\qquad(n\ge1).
\]
The OEIS entry~\cite{OEIS} states that $a(n)$ is even for $n\le400$; this is
not used in the proofs. For general Gaussian parts with at least two distinct split primes,
the component vectors are no longer collinear, and the rank-one Farey
argument of Section~\ref{sec:parity} alone does not apply.

Theorem~\ref{thm:power-growth} gives
$\limsup_{e\to\infty}a(N^e)/e\ge2\sqrt2$, while
Theorem~\ref{thm:polynomial-bound} provides the general polynomial
upper bound $a(N^e)=O_N(e^{\omega(N)})$.
For prime powers $q^e$ with $q\equiv1\pmod4$, the sharp growth
order is $\Theta(e)$. The existence of a limit
$\lim_{e\to\infty}a(q^e)/e$, its value if it exists, and
precise asymptotics for other fixed bases remain open.
In particular, the available upper and lower estimates do not determine
the order of growth of $a(2^e)$. It is also open how the factorization
of $n$ determines $a(n)$ beyond the semiprime case.

\enlargethispage{5\baselineskip}
\section*{Acknowledgments}
The author thanks the anonymous referee for a careful and constructive
report, including suggestions leading to the general polynomial bound,
its prime-power specialization, the strengthened growth constant,
and the rigidity of exceptional extremisers.

\section*{Statements and Declarations}

\paragraph{Funding.}
The author declares that no funds, grants, or other support were received
for the preparation of this manuscript.

\paragraph{Competing interests.}
The author has no relevant financial or non-financial interests to disclose.

\paragraph{Author contributions.}
Felix Huber conceived and directed the study and selected the mathematical
arguments and conclusions. He accepts full responsibility for the accuracy
and integrity of the manuscript, including its mathematical statements,
proofs, references, and figures, and for the final submission.

\paragraph{Data and code availability.}
No external datasets, software, or numerical certificates are required for
the results. All proofs and constructions are contained in the manuscript.
OEIS A398489 and A398490 are cited for context only.

\paragraph{Use of artificial intelligence.}
OpenAI ChatGPT, Anthropic Claude, and Microsoft Copilot were used during
August--November 2026 for proof exploration, independent critical assessments,
numerical cross-checks, drafting, revision, and language editing.
In the November revision, ChatGPT assisted in incorporating the referee's
suggestions and preparing mathematical and editorial amendments;
Claude independently audited the revised arguments and identified further
clarifications. AI assistance also supported TikZ figure preparation.
Numerical checks are not used as proof. The author retains responsibility
for independently checking and approving all mathematical statements,
proofs, calculations, references, figures, and declarations before submission.

\end{document}